\documentclass{amsart}
\usepackage{graphicx}
\usepackage{amsthm,amsfonts}
\usepackage{amsmath}
\usepackage[latin1]{inputenc}
\usepackage{hyperref}
\usepackage{amssymb}
\usepackage{latexsym}
\usepackage{psfrag}

\newtheorem{thm}{Theorem}[section]

\newtheorem{lem}[thm]{Lemma}
\newtheorem{prop}[thm]{Proposition}

\theoremstyle{definition}
\newtheorem{remark}[thm]{Remark}

\newcommand{\disp}{\displaystyle}

\def \build#1#2#3{\mathrel{\mathop{\kern 0pt#1}\limits_{#2}^{#3}}}

\newcommand{\ac}[1]{\left\{#1\right\}}
\newcommand{\pa}[1]{\left(#1\right)}
\newcommand{\cro}[1]{\left[#1\right]}
\newcommand{\abs}[1]{\left|#1\right|}
\newcommand{\norm}[1]{\left\Vert#1\right\Vert}

\newcommand{\esp}[1]{\mathbb{E}\left[#1\right]}
\newcommand{\pr}[1]{\mathbb{P}\left(#1\right)}
\newcommand{\var}[1]{\textnormal{Var}\left(#1\right)}

\newcommand{\egenloi}{\,\build{=}{}{\text{law}}\,}

\newcommand{\gl}{\lambda}

\newcommand{\bigO}[1]{\mathcal O\pa{#1}}
\newcommand{\To}{\longrightarrow}

\newcommand{\esup}[1]{\lceil#1\rceil}
\newcommand{\scrip}{\scriptstyle}
\newcommand{\lr}{\mathcal{L}}

\newdimen\AAdi%
\newbox\AAbo%
\def\AArm{\fam0 \rm}%
\def\AAk#1#2{\setbox\AAbo=\hbox{#2}\AAdi=\wd\AAbo\kern#1\AAdi{}}%
\def\AAr#1#2#3{\setbox\AAbo=\hbox{#2}\AAdi=\ht\AAbo\raise#1\AAdi\hbox{#3}}%
\def\BBone{{\AArm 1\AAk{-.8}{I}I}}%

\newcommand\ie{i.e.\spacefactor=1000}
\newcommand\set[1]{\ensuremath{\{#1\}}}
\newcommand\xfrac[2]{#1/#2}
\newcommand\bigpar[1]{\bigl(#1\bigr)}

\newcommand\E{\mathbb{E}}
\newcommand\Bi{\operatorname{Bi}}

\newcommand\sumkj{\sum_{k=1}^\infty\sum_{j=1}^{2^k}}
\newcommand\ikj{I_{k,j}} 
\newcommand{\gd}{\delta}

\title{Quicksort with Unreliable Comparisons: A Probabilistic Analysis}

\author[Alonso, Chassaing, Gillet, Janson, Reingold \& Schott]{Laurent Alonso
\and Philippe Chassaing
\and Florent Gillet
\and Svante Janson
\and Edward M. Reingold
\and Ren\'e Schott
}

\address{INRIA-Lorraine and LORIA,
Universit\'e Henri Poincar\'e-Nancy I,
BP 239, 54506,
Van\-d{\oe}uvre-l\`es-Nancy, France}
\email{\tt Laurent.Alonso@loria.fr}

\address{Institut \'{E}lie Cartan,
Universit\'e Henri Poincar\'e-Nancy I,
BP 239, 54506,
Vand{\oe}uvre-l\`es-Nancy, France}
\email{Philippe.Chassaing@antares.iecn.u-nancy.fr}
\urladdr{http://www.iecn.u-nancy.fr/\~{}chassain/}

\address{Institut \'{E}lie Cartan,
Universit\'e Henri Poincar\'e-Nancy I,
BP 239, 54506,
Vand{\oe}uvre-l\`es-Nancy, France}
\email{Florent.Gillet@antares.iecn.u-nancy.fr}

\address{Department of Mathematics,
Uppsala University,
PO Box 480,
S-751 06 Uppsala, Sweden}
\email{svante.janson@math.uu.se}
\urladdr{http://www.math.uu.se/\~{}svante/}

\address{Department of Computer Science,
Illinois Institute of Technology,
Stuart Building,
10 West 31st Street, Suite 236,
Chicago, IL  60616-3729  U.S.A.}
\email{reingold@iit.edu}

\address{Institut \'{E}lie Cartan and LORIA,
Universit\'e Henri Poincar\'e-Nancy I,
BP 239, 54506,
Vand{\oe}uvre-l\`es-Nancy, France}
\email{schott@loria.fr}

\begin{document}

\begin{abstract}
We provide a probabilistic analysis of the output of Quicksort
when comparisons  can err.
\end{abstract}


\maketitle

\section{Introduction}

Suppose that a sorting algorithm, knowingly or unknowingly,
 uses element comparisons that can
err. 
Considering sorting algorithms based solely on binary
comparisons of the elements to be sorted (algorithms such as insertion sort,
selection sort, quicksort, and so on), what problems do we face when those
comparisons are unreliable?  For example, \cite{spot-checkers} gives a clever
$\bigO{\epsilon^{-1} \log n}$
algorithm to assure, with probability $1-\epsilon$,
that a putatively sorted sequence of length $n$ is truly sorted.  But  knowing
the structure of the ill-sorted output  would likely make error
checking easier. Also,  in situations
in which a reliable comparison is the fruit of a long process, 
one could chose to interupt the comparison process, 
thus trading reliability of comparisons (and quality of the output) for time.
As a first step in order to understand the consequences of errors,  
we propose to analyze
the number  of inversions in the output
of a sorting algorithm (we choose Quicksort  \cite{Knuth3}) subject to errors.

We assume throughout this paper that the elements of the sequence
\[x=(x_1,x_2,\dots,x_n)\]
   to  be sorted are distinct.
We assume further that the only comparisons subject to error are those made
between elements being sorted; that is, comparisons among indices and so on
are always correct.  Errors in element comparisons are random events,
spontaneous and
independent of each other, of position, and of value, with a common
probability
   $p$, $n$   being the length
   of the list to be sorted. 
   The number of inversions
   in the output sequence
$y=(y_1,y_2,\dots,y_n)$ is denoted
\[I(y)=\#\ac{(i,j)\ \left\vert\ 1\le i<j\le n\mbox{ and }y_i>y_j\right.}.\]
We assume that the input list is presented in random order,
each of the $n!$ random orders being equiprobable.
Finally we denote by $I(n,p)$ the random number of inversions
in the output sequence of Quicksort subject to errors.

Our result is, roughly speaking,
\[I(n,p)=\Theta\left(n^2p\right),\]
when $(n,p)\to(\infty,c)$,  meaning that
$\frac{I(n,p)}{n^2p}$  converges  to some nondegenerate
probability distribution.
The ``surprise", not so unexpected after the fact,
   is that there are  phase changes in the limit law, depending on the
asymptotic behaviour of
$(n,p)$.

The organization of this paper is as follows:
The results are stated in Section \ref{S:results}.
In Section \ref{functequns},
we  establish a
general distributional identity
for $I(n,p)$. In the remaining sections,
  we prove convergence results for $I(n,p)$ when:
\begin{itemize}
\item $p\rightarrow c$,  $0<c\le 1$,
\item $p$ vanishes more slowly than $1/n$,
\item $p\sim \gl/n$ where  $\gl$ is a positive constant.
\end{itemize}
The case $np\to0$ is different and not treated in detail; see
Remark \ref{R:toosmall}.
In Section \ref{fixedpts}, we establish  a general result of convergence
using contraction methods
(cf.\ \cite{Rosler91,RoeRue}), and we use it in Section
\ref{proofTh1},
for the first two cases.
These methods do not  apply  for Case $3$, which requires
poissonization (see Section \ref{proofTh3}, where
we use an embedding of Quicksort in a Poisson point process).

\section{Results}
\label{S:results}

Set
\[X_{n,p}=\frac{I(n,p)}{n^2p}.\]
We will always let $U$
denote a random variable that is uniformly
distributed on $[0,1]$. Also,  ${\mathbb N}^*$
shall denote the set of positive integers, and ${\mathbb N}$
  the set of nonnegative integers.

\subsection{Case $1$: $\lim p= c> 0$}
\begin{thm}
\label{p=c}
If\/ $\lim p=c$, $c\in(0,1]$, then
$X_{n,p}$ converges in distribution
to a random variable
$X_c$ whose distribution is characterized as the unique
solution with finite mean of the equation
\begin{align}
\label{p=Ode1}
X_c\overset{\text{law}}{=}
[(1-2c)U+c]^2X_c+[(2c-1)U+1-c]^2\widetilde{X}_c+T(c,U),
\end{align}
in which
$\widetilde X_c$ denotes a copy of $X_c$,
$(X_c,\widetilde X_c, U)$  are independent, 
and
\begin{align*}
T(c,U)=\frac{1-c}{2}(U^2+(1-U)^2)+cU(1-U).
\end{align*}
Furthermore,
\[
\esp{X_c}= \frac{2-c}{2(1+2c-2c^2)},\]
and
\[\var{X_c}=
\frac{(1-c)^2(1-2c)^2}{4(1+2c-2c^2)^2(3+6c-8c^2+4c^3-2c^4)}.
\]
\end{thm}

As usual with laws
related to Quicksort, see e.g.\
\cite{Rosler91,RoeRue},
$nU$ is approximately the position
of the pivot of the first step of the algorithm.
As in standard Quicksort recurrences,
the coefficients of $X_c$ and of its independent copy
$\tilde X_c$ are related to the sizes of the two sublists
on the left and right of the pivot,
sizes respectively asymptotic
to $n\pa{(1-2c)U+c}$
and $n\pa{(2c-1)U+1-c}$.
The toll function
$T(c,U)$ is approximately $(n^2p)^{-1}\approx(n^2c)^{-1}$
times the number of inversions created in the first step:
$\xfrac{c(1-c)n^2U^2}2$ is approximately
the number of inversions of the $cnU$
elements, smaller than the pivot but misplaced
on the right of it, with the $\left(1-c\right)nU$ elements
 smaller than the pivot, that are placed,
as they should be,
on the left;   ${c^2n^2U(1-U)}$ is the number of inversions
between misplaced elements from
the two sides of the pivot.
The toll function $T(c,U)$ depends on only one of the two sources of randomness
(the randomly ordered input list, and the places of the errors),
viz.,\ the first one, through $U$. The second source of randomness
is killed by the law of large numbers:
in the average, each of the $cnU+o(n)$ misplaced
numbers from the right of the pivot
produces inversions with one half of the $(1-c)nU+o(n)$
elements smaller than the pivot, that are placed,
as they should be,
on the left. As opposed to the other values of $c$,
the choices $c=0.5$
and $c=1$ lead to deterministic $X_c=1/2$,
without any surprise : for $p=0.5$  the output sequence is
   a random uniform permutation, with a number of inversions
concentrated around
$n^2/4$ \cite[Chap. 5.1.1]{Knuth3}; for $p=1$  the output sequence is
   decreasing, and has $n(n-1)/2$ inversions.

\subsection{Case $2$:
$p$ vanishes more slowly than $\frac1n$}
\begin{figure}
\begin{center}
\includegraphics[width=4cm]{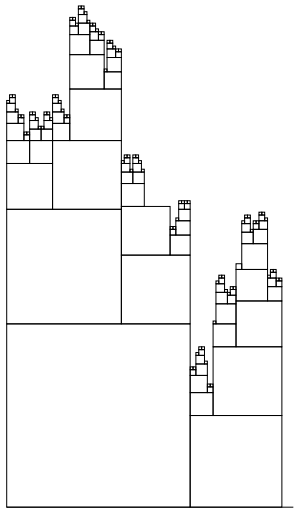}
\end{center}
\caption{The Find process.}
\label{Findprocess}
\end{figure}

\begin{thm}
\label{npinfini}
If $\lim p= 0$ and $\lim np=  +\infty$,
$X_{n,p}$ converges in distribution  to a random variable
$X$ whose distribution is characterized as the unique
solution with finite mean of the equation
\begin{align}
\label{p=ode1}
X&\overset{\text{law}}{=}U^2X+(1-U)^2\widetilde{X}+
\frac{U^2+(1-U)^2}{2}.
\end{align}
In \eqref{p=ode1}, $\widetilde X$ denotes a  copy of $X$
and $(X,\widetilde X,U)$  are independent.
Furthermore,
\begin{align*}
\esp{X}=1\qquad\text{and}\qquad
\var{X}=\frac{1}{12}.
\end{align*}
\end{thm}
Note that equation (\ref{p=ode1})
is just     (\ref{p=Ode1})
specialized to  $c=0$, but, as opposed to $c\neq 0$,
an additional condition,   $p\gg 1/n$,
is needed to ensure that
 the law of large numbers  still holds.
Also, as another difference 
between (\ref{p=Ode1}) and (\ref{p=ode1}), for $p\ll 1$
   the errors do not change  the sizes of the sublists
 in a significant way.
The solution $X$ equals half  the sum of the squares of
the widths of the random  intervals $[Y_{k,j},Y_{k,j+1}]$
  defined
  by (\ref{recurrence_sur_U})  below.
This is equivalent to the following statement:
\begin{prop}
\label{propfindprocess}
The solution $X$ equals half  the area $\int_0^1 Z(t)\,dt$ under the
FIND limit
process $Z$.
\end{prop}
\noindent For both these claims,
see Remark \ref{R:1}. The Find process was introduced in \cite{Gruber_Rosler96}
and is pictured at Figure \ref{Findprocess}.

\subsection{Case $3$: $\lim np={\gl}$}
\quad

Assume that
\begin{itemize}
\item
$\Pi$ is a Poisson point process with intensity $\gl$ on
${\mathbb N}^*\times[0,1]$,
meaning that, for each $n$,  $\vert\Pi\cap\pa{\{n\}\times[0,1]}\vert$  is a Poisson
random variable with mean $\gl$, and the second coordinates of  points of $\Pi$
are uniform on $[0,1]$ and independent (see
\cite{PPP} for a  general definition of Poisson point processes);
\item
$\{U_{k,j}: k\ge 0, 1\le j\le 2^k\}$
is an array of independent  uniform random variables on $[0,1]$,
independent of $\Pi$;
\item
the random variables  $(Y_{k,j},k\geq 0,1\leq j\leq 2^k)$ are
defined recursively by
\begin{align}\label{recurrence_sur_U}
\begin{array}{c}
Y_{0,0}=0,\qquad Y_{0,1}=1,\qquad Y_{k+1,2j}=Y_{k,j}\qquad
\mbox{for } 0\leq j\leq2^k,\\
Y_{k+1,2j-1}=(1-U_{k,j})Y_{k,j-1}+U_{k,j}Y_{k,j}\qquad
\mbox{for }1\leq j\leq 2^k;
\end{array}
\end{align}
\item
for $x\in[0,1]$, $J_k(x)=2j-1$ if $Y_{k-1,j-1}\leq x < Y_{k-1,j}$,
\end{itemize}
and define, for $\gl>0$,
(the sum is a.s.\ finite by Lemma \ref{L:mean})
\begin{align}
\label{def de I}
X(\gl)
=
\frac1{\gl}\ \sum_{(k,x)\in\Pi}\big|x-Y_{k,J_k(x)}\big|.
\end{align}

The variables $Y_{k,j}$ describe a fragmentation process
(see
 \cite{Gruber_Rosler96} for historical references):
we start with $[0,1)$ and recursively break each interval into two at
a random point (uniformly chosen). In the $k$-th generation we thus
have a partition of $[0,1)$ into $2^k$ intervals $\ikj$, $1\le j \le
2^k$,
with $\ikj=[Y_{k,j-1},Y_{k,j})$.
The interval of generation $k-1$ that
contains $x$ is cut at step $k$ at the point $Y_{k,J_k(x)}$.
Hence $\big|x-Y_{k,J_k(x)}\big|$ in \eqref{def de I} is the distance
from $x$ to this cut point.

\begin{thm}\label{p=1/n}
If\/ $\lim p= 0$ and\/ $\lim np=  \gl>0$,
then
$X_{n,p}$ converges in distribution  to
$X(\gl)$.
The family $\set{X(\gl)}_{\gl>0}$
of random variables satisfies
the distributional identity:
\begin{equation}
\label{functional1}
X(\gl)\build{=}{}{\text{law}}U^2 X(\gl U) + (1-U)^2  \widetilde X(\gl(1-U))
+ \Theta(\gl,U),
\end{equation}
in which, conditionally given that $U=u$,
$X(\gl U)$, $\widetilde X(\gl(1-U))$ and  $\Theta(\gl,U)$
are independent,
$X(\gl U)$ and $\widetilde X(\gl(1-U))$
are distributed  as  $X(\gl u)$ and $X(\gl(1-u))$,
respectively, and
\begin{align*}
\Theta(\gl,u)\build{=}{}{\text{law}}\frac1\gl\sum_{i=1}^{N_\gl}\abs{u-V_i},
\end{align*}
in which  $N_\gl$ is a Poisson random variable with mean $\gl$,
the random variables
$V_i$ are uniformly distributed on $\left[0,1\right]$, and
$N_\gl$, and the  $V_i$'s are independent.
Furthermore,
\begin{align}
\esp{X(\gl)}=1, \qquad \var{X(\gl)}=\frac1{12}+\frac1{3\gl}.
\end{align}
\end{thm}

\begin{remark}
\label{remmesure}
The distributional identities \eqref{p=Ode1}, \eqref{p=ode1} and
\eqref{functional1} really are equations for distributions,
but it is more convenient to state them for random variables as done
here.
For \eqref{functional1} to make sense, i.e. in order to insure, for instance,  that
$X(\gl U)$ is  a random variable,  it is implicitly assumed that
the random variables $X(\gl)$ depend measurably on $\gl$.
Thus a solution of \eqref{functional1}  is a family of probability measures
$\mu=\pa{\mu_\gl}_{\gl>0}$ on $[0,+\infty)$, such that there exists 
 a family $Y=\pa{Y(\gl)}_{\gl>0}$ 
of random variables defined on the same probability space
$(\Omega,\mathcal A,\mathbb P)$ 
and satisfying the following properties:
\begin{enumerate}
\item  for $\gl>0$, $\mu_\gl$ is the distribution of $Y(\gl)$,
\item  $Y$  is a measurable process \cite[Chap. 1]{Karatzas}, meaning that
the mapping
\[
(\gl,\omega)
\rightarrow
Y(\gl,\omega)
:\ \ 
\pa{(0,+\infty)\times\Omega,\mathcal B((0,+\infty))\otimes\mathcal A}
\rightarrow
\pa{[0,+\infty),\mathcal B([0,+\infty))}
\]
is measurable,
\end{enumerate}
and such that \eqref{functional1} holds for $Y$.
A measurable version of the stochastic process
$X=\pa{X(\gl)}_{\gl>0}$ 
is defined at
\eqref{samuel} below (measurability follows  from 
\cite[Rem. 1.14]{Karatzas}).
\end{remark}

For uniqueness, we need extra assumptions:
let $\mathcal M$
denote the class of families of distributions $\mu=\pa{\mu_\gl}_{\gl>0}$
satisfying (i) and (ii) above,  plus the condition:
\begin{enumerate}
\item[(iii)]  for some $\alpha\in\left(0,1\right)$,  the function
$$
\gl
\longrightarrow
\lambda^\alpha \esp{Y(\gl)}
$$
is bounded on any bounded interval  of $(0,+\infty)$.
\end{enumerate}
Let $\nu_\gl$ denote the distribution of $X(\gl)$. We have

\begin{thm}
\label{thm:unique} 
The family $\nu=\pa{\nu_\gl}_{\gl>0}$
is the unique solution of \eqref{functional1}
in $\mathcal M$.
\end{thm}

We do not know whether the extra assumption (iii) is necessary.
Let us comment further on equation \eqref{def de I}.
Writing $\Pi_k=\set{x:(k,x)\in\Pi}$ and $\Pi_{k,j}=\Pi_k\cap\ikj$, we
can thus rewrite \eqref{def de I} as
\begin{equation}
\label{emma}
X(\gl)=\frac1{\gl}\sumkj \sum_{x\in\Pi_{k,j}}|x-x_{k,j}|,
\end{equation}
where $x_{k,j}$ is either the left or right endpoint of $\ikj$
(depending on whether $j$ is even or odd).

Note that, conditioned on the partitions $\set{\ikj}$, \ie{}
on $\set{Y_{k,j}}_{k,j}$, each $\Pi_{k,j}$ is a Poisson process on
$\ikj$ with intensity $\gl$, with the processes $\Pi_{k,j}$ independent.
Since only the distribution of $X(\gl)$ matters, we can by this
conditioning and an obvious symmetry of the Poisson processes $\Pi_{k,j}$
just as well let $x_{k,j}$ in \eqref{emma}
be the left endpoint of
$\ikj$ for every $k$ and $j$.

Let $\Pi'$ be a Poisson process on $(0,1]\times(0,\infty)$
with intensity $1$,
and let $\xi(t)=\sum_{(x,y)\in\Pi',\,y\le t} x$, $t\ge0$.
(This is a pure jump L\'evy process with L\'evy measure $\BBone_{(0,1]} dt$.)
Let $\xi^{(k,j)}(t)$ be independent copies of this process,
independent of $\set{Y_{k,j}}$.
A scaling argument shows that \eqref{emma} can be written
\begin{equation}
\label{samuel}
X(\gl)=\frac1{\gl}\sumkj |\ikj| \xi^{(k,j)}(\gl|\ikj|).
\end{equation}

\begin{remark}\label{R:1}
Let $\hat X=\tfrac12\sumkj|\ikj|^2$. Then $\hat X$ satisfies \eqref{p=ode1},
so that  $\hat X$  is the limit variable $X$ in Theorem \ref{npinfini}.
($X$ is a.s.\  finite and has finite mean by Lemma \ref{L:erika}.)
Moreover,
the FIND limit
process $Z$ in \cite{Gruber_Rosler96} is defined by
$Z(t)=\sumkj|\ikj|\BBone_{t\in\ikj}$; hence
$\int_0^1 Z(t)\,dt=\sumkj|\ikj|^2=2X$.
This justifies Proposition \ref{propfindprocess}.

Moreover, by the law of large numbers,
$\E|\gl^{-1}\xi(\gl)-1/2|\to0$
as $\gl\to\infty$. It follows  (by dominated
convergence using Lemma \ref{L:erika})
that for the special version of $X(\gl)$ 
defined at \eqref{samuel} 
\[\E|X(\gl)-X|\to0,\]
and hence
$X(\gl)$ converges to $X$ in distribution as $\gl\to\infty$.
\end{remark}

In this third case, we have a system of equations involving
an infinite  family of laws, and we could not adapt
  the contraction method: we rather use a poissonization.
  The phase transition from (\ref{p=ode1}) to (\ref{functional1})
is explained easily: instead of a number of errors $\gg 1$,
we have now  $\bigO{1}$  errors at each step, and the  law of large numbers
does not hold anymore for the number of inversions produced by step 1. Actually
the number $N_\gl$ of errors at the first step  is
  asymptotically  Poisson distributed, and the $N_\gl$ errors are
at positions  $nV_i$, approximately uniformly
distributed on $[0,1]$.  Thus, the number of inversions
caused by this first step is approximately
\[n \sum_{i=1}^{N_\gl} \vert U-V_i\vert
\approx n^2p \Theta(\gl,U)
.\]

\begin{remark}
Actually we prove a stronger theorem in each of the three cases, as
we prove convergence of laws for the Wasserstein $d_1$ metric \cite{Rachev}.
It entails convergence of the first moment.
The convergence of higher moments is an open problem.
\end{remark}

\begin{remark}
As we shall see in Section \ref{proofTh3},
the distribution tail $\pr{X(\gl)\geq x}$ decreases exponentially fast
  (Theorem \ref{tail}).
\end{remark}

\begin{remark}
\label{R:toosmall}
When $np\to0$ very slowly, that is $(np)^{-1} \ll \log n$, we conjecture that
$2np\log \pa{I(n,p)/n}$ converges in distribution to $\log U$,
  with the consequence that $n^{1-\varepsilon }\ll I(n,p)\ll n$,
for any positive $\varepsilon$. Actually, the main contribution to  $I(n,p)$
comes from the "first" error, in some sense. When
$(np)^{-1} \sim \log n$, the probability
  that no error occurs has a positive limit: we conjecture that, conditionally
given the occurence of at least one  error, the situation is similar
  to the previous case, that
is, $\log \pa{I(n,p)}/\log n$ converges in distribution to
   a random variable  with values in  $(0,1)$. When $(np)^{-1} \gg \log n$,
$\pr{I(n,p)=0}\to 1$.
\end{remark}

\begin{remark}
\label{R:filter}
Finally,  we would like to stress that
in the proof of convergence for one the three regimes  considered
in this Section, we have to deal 
simultaneously
with any sequence $(n,p_n)$ converging to $(+\infty,c)$ 
according to this regime. This can be observed on the key equation
\eqref{eq en I}, for instance, in which  we would like to argue,
roughly speaking, that if 
$(n,p)$ is close to  $(+\infty,c)$ 
according to  a given regime, then  $(Z_{n,p}-1,p)$
and $(n-Z_{n,p},p)$ are also close to  $(+\infty,c)$ 
according to the same  regime, with a large probability: here the same probability
$p$ is associated to three different integers,
 $n$, $Z_{n,p}-1$
and $n-Z_{n,p}$, that denote the sizes of the input list, 
and of the two sublists formed at the first step of Quicksort,  respectively.
Thus $p$ cannot be seen as a sequence indexed by $n$.
In order to allow such a loose relation
between $n$ and $p$, filters turn out to be more handy than sequences
  (see \cite[Chap. I]{Bourbaki}). Convergences in the three regimes are
 thus understood as convergences along the three corresponding filters
(see Theorem \ref{th de convergence}).
\end{remark}

\section{A distributional identity for the  number of inversions}
\label{functequns}
At the first step Quicksort compares all elements of the input list with the
   first element of the list (usually called \emph{pivot}). All items
less (resp.\ larger)
   than the pivot are stored in
a sublist on the left (resp.\ right) of the pivot.
Comparisons are not reliable, therefore $s_{\ell}$ items that should
belong to the
  left sublist are wrongly
stored in the right sublist, and
$s_{r}$ items larger than the pivot are  misplaced in the left sublist.

Since its  items are chosen randomly, the input list is a random permutation
and the true rank of the pivot    can be
written $\esup{nU}$, where $U$ is uniformly distributed on $[0,1]$ and
$\esup{x}$ is the ceiling of  $x$. Also, conditionally given $U$,
$s_{\ell}$ (resp.\ $s_{r}$)  is a binomial random variable with
parameters
($\esup{nU}-1$, $p$)
(resp.\ ($n-\esup{nU}$, $p$)).
Quicksort with error is then independently applied on the left sublist
$\ell$ and on the right sublist  $r$
and new errors occur, ultimately producing
two new sublists $\tilde{\ell}$ and  $\tilde{r}$.
Set
\begin{equation*}
Z_{n,p}=\lceil nU \rceil-s_{\ell}+s_{r},
\end{equation*}
so that $Z_{n,p}-1$ (resp.\ $n-Z_{n,p}$)  is the
size of $\ell$ and $\tilde {\ell}$  (resp.\ $r$ and  $\tilde{r}$).

In order to enumerate the inversions of the output list,
we introduce a \textit{purely fictitious} error-correcting
   algorithm
that parallels the implementation  of Quicksort:
This fictitious error-correcting
   algorithm  has two recursive steps,

\begin{itemize}
\item
First, the error-correcting
   algorithm corrects the sublists $\tilde  \ell$ (resp.\  $\tilde{r}$)
   at costs $L=I(\tilde \ell)$  (resp.\ $R=I(\tilde{r} $)), producing two
increasing
sublists $\hat  \ell$ and $\hat{r}$. Note that $L$ and $R$ are conditionally
   independent, given $Z_{n,p}$. Furthermore, the two sublists $\ell$
and $r$ obtained at the end of Step $1$ are  in
   uniform random order before the second step of Quicksort,
so that, conditionally given $Z_{n,p}$, cost $L$
(resp.\ $R$) is distributed as $I(Z_{n,p}-1,p)$ (resp.\ $I(n-Z_{n,p},p)$).

\item  Then the error-correcting
   algorithm corrects the errors of Step $1$, at a cost
   $t(n,p)=I(\hat{\ell}\Vert\textrm{pivot}\Vert\hat r)$.
Here $\hat{\ell}\Vert\textrm{pivot}\Vert\hat r$ stands for the list obtained
when one puts $\hat{\ell}$, the pivot and $\hat r$
side by side. The number of
inversions $t(n,p)$
in the list
$\hat{\ell}\Vert\textrm{pivot}\Vert\hat r$ is analyzed in detail at the end
of this section.
\end{itemize}

\begin{figure}[ht]
\label{2sublists0}
\begin{center}
\psfrag{sl}{$s_{\ell}$}
\psfrag{sr}{$s_r$}
\psfrag{r}{$r$}
\psfrag{l}{$\ell$}
\psfrag{a}{$\widetilde{\ell}$}
\psfrag{b}{$\widehat{\ell}$}
\psfrag{c}{$\vec{\ell}$}
\psfrag{e}{$\widetilde r$}
\psfrag{f}{$\widehat r$}
\psfrag{g}{$\vec r$}
\psfrag{i}{$W_{\vec{\ell}}$}
\psfrag{j}{$W_{\vec r}$}
\psfrag{k}{$s_{\ell} s_r+s_{\ell}+s_r$}
\psfrag{Izn}{$I(Z_{n,p}-1,p)$}
\psfrag{Inzn}{$I(n-Z_{n,p},p)$}
\psfrag{Quicksort}{Quicksort}
\psfrag{T(n,p)}{$T(n,p)$}
\psfrag{I(n,p)}{$I(n,p)$}
\psfrag{h}{first step}
\psfrag{Algorithme de}{correction}
\psfrag{correction}{algorithm}
\psfrag{w}{$\lceil nU\rceil$}
\includegraphics[width=12cm]{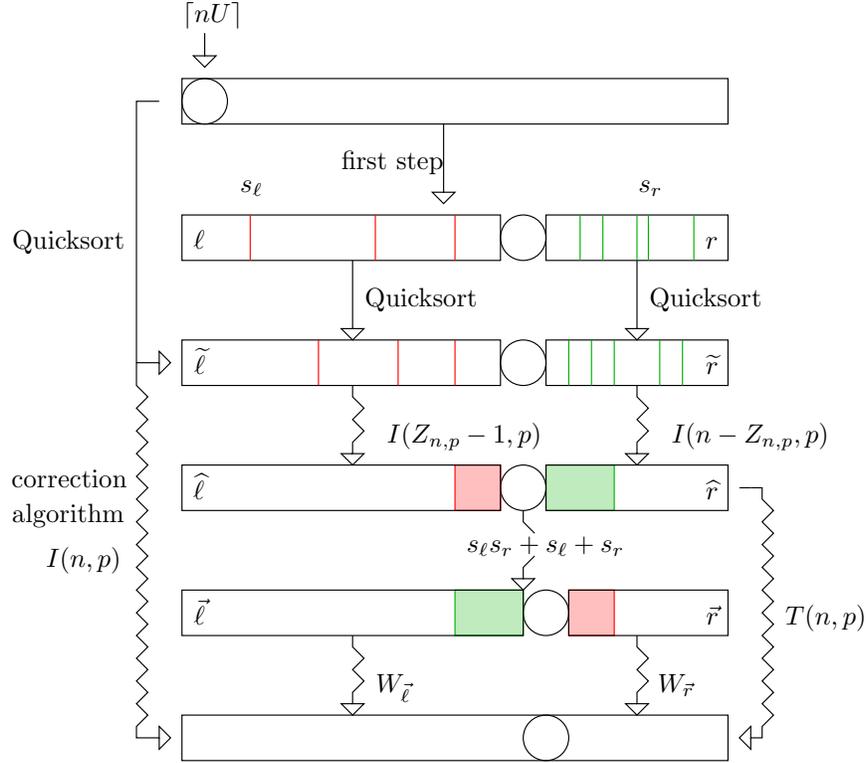}\\
\end{center}
\caption{The error-correcting
   algorithm.}
\end{figure}

\smallskip

These two steps lead to the following equation for $I(n,p)$:
\begin{equation}\label{eq en I}
I(n,p)\build ={}{law}I(Z_{n,p}-1,p)+
I'(n-Z_{n,p},p)+
t(n,p)
\end{equation}
where $Z_{n,p}=\lceil nU \rceil-s_{\ell}+s_{r}$.  We shall obtain
the asymptotic distribution of $t(n,p)$, and as a consequence
(\ref{eq en I})  will translate, after renormalisation,
   into a distributional identity satisfied by the limit law of
$I(n,p)/(n^2p)$. The limit law  appears  on both sides of   the
distributional identity,
as expected, due to the recursive structure of Quicksort, and is thus
  characterized as the fixed point of some transformation.

\subsection{Description of $t(n,p)$}

\begin{figure}[ht]
\label{2sublists}
\begin{center}
\includegraphics[width=12cm]{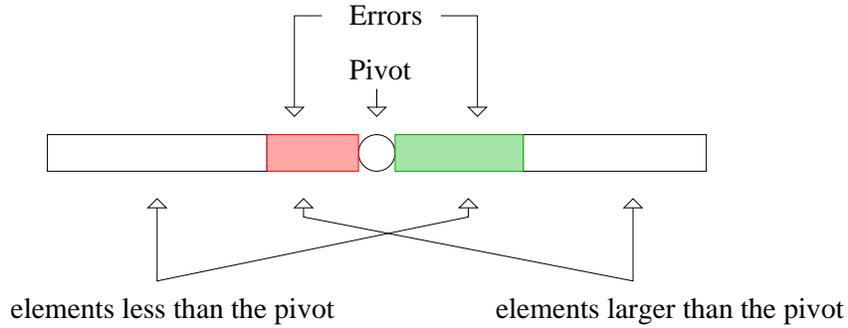}
\end{center}
\caption{The two sublists $\hat \ell$ and $\hat r$.}
\end{figure}

At the end of the first step of the error-correcting
   algorithm, we obtain two
subarrays  $\hat \ell$ and $\hat r$,
left and right of the pivot  (cf.\ Figure \ref{2sublists}).
They are sorted in increasing order but there are $s_{r}$ (red)
elements larger than
the pivot just to its left and $s_{\ell}$ (green) elements smaller than
the pivot element
just to its right.  Thus, the only misplaced elements that the
proofreader must correct in
step 2 are clustered around the pivot.

In order to sort the list, the red and green sublists must be
exchanged. This requires
$s_{\ell}s_{r}+s_{\ell}+s_{r}$ inversions. We get therefore two unsorted lists
$\vec\ell$ and $\vec r$
each composed of two sorted
sublists. All items of $\vec{\ell}$  (resp.\ of $\vec r$ )  are now
smaller (resp.\ larger)
than the pivot, so that  the  length of $\vec\ell$  (resp.\ of $\vec r$ )  is
$\esup{nU}-1$ (resp.\ $n-\esup{nU}$).
It remains to sort $\vec\ell$ and $\vec r$, at respective costs
$W_{\vec\ell}$  and
${W}_{\vec r}$ that are conditionally independent given $U$,
leading to:
\begin{equation}
\label{t1np}
t(n,p)=s_{\ell}s_{r}+s_{\ell}+s_{r}+W_{\vec\ell}+{W}_{\vec r}.
\end{equation}

\subsection{A model for $\pa{W_{\vec\ell},{W}_{\vec r}}$}
Let $W_m$ be the number of inversions  in a list of $m$ elements
   sorted as follows: each element is
painted black (white) with probability $p$ (resp.\ $1-p$). Then the
black and white sublists
are separately sorted in increasing order
and the two sorted  sublists are placed side by side, producing a
new list $h$
with $m$ elements.  We have
\begin{prop}
\label{invb&w}
  Let $Y_1,\dots,Y_m$ be $m$ independent
   Bernoulli random variables with the same parameter
$p$, and let $S_m=Y_1+\ldots+Y_m$. Then
\[
W_m\build ={}{law}\Big(\sum_{i=1}^miY_i\Big)-\frac{S_m(S_m+1)}{2}.
\]
\end{prop}

\begin{proof}
Let us abbreviate $S_m$ to $S$. Among the $Y_i$'s, let
$Y_{i_1},\dots,Y_{i_S}$ denote
   the $S$ random variables equal to $1$,
$ Y_{i_{S+1}},\ldots,Y_{i_m}$ those equal to $0$, with $i_1<\ldots<i_S$ and
$i_{S+1}<\ldots<i_m$. Now $W_m$ can be seen as the number of inversions of
the list $\big(i_j\big)_{1\leq j\leq m}$.
In order to move the numbers $i_j$ ($j\leq S$) to the correct
position, the proofreader
corrects  inversions
with each of the $i_j-j$  elements of $\{1,\dots,m\}$ that are
smaller than  $i_j$ and  do
not belong to $\{i_1,\ldots,i_S\}$. Thus
\begin{equation}
\label{poissonstep1}
W_m=\sum_{j=1}^{S}(i_j-j),
\end{equation}
leading to the result.
\end{proof}

With the help of Proposition \ref{invb&w}, we can give a useful description
  of the distribution of $(s_{\ell},W_{\vec\ell})$ and   $(s_{r},W_{\vec r})$:
\begin{prop}
\label{description}
Conditionally given that the length of
$\vec{\ell}$  is $m-1$, $(s_{\ell},W_{\vec\ell})$ and
$(s_{r},W_{\vec r})$
  are independent and distributed as $(S_{m-1},W_{m-1})$ and
  $(S_{n-m},W_{n-m})$, respectively.
\end{prop}

To sum up the results of this section, renormalizing (\ref{eq en I}),
one obtains
a distributional identity satisfied by $\disp X_{n,p}$:
\begin{equation}\label{eq en X}
X_{n,p}\build ={}{law}A_{n,p} X_{Z_{n,p}-1,p}+B_{n,p}
\widetilde{X}_{n-Z_{n,p},p}+T_{n,p}
\end{equation}
in which
\begin{align}
Z_{n,p}&=\lceil nU \rceil-s_{\ell}+s_{r}, \label{znp}
\\
A_{n,p}&=\pa{\frac{Z_{n,p}-1}{n}}^2,
\\
B_{n,p}&=\pa{\frac{n-Z_{n,p}}{n}}^2,
\\
t(n,p)&=s_{\ell}s_{r}+s_{\ell}+s_{r}+W_{\vec\ell}+W_{\vec r} ,
\\
T_{n,p}&=\frac{t(n,p)}{n^2p}, \label{tnp}
\end{align}
and
\begin{itemize}
\item $U$ is a uniform random variable on $[0,1]$, and
$\esup{nU}$ is the position of the pivot,
\item conditionally given $\esup{nU}=m$,
$(s_{\ell},W_{\vec\ell})$ and   $(s_{r},W_{\vec r})$
are distributed as in Proposition \ref{description},
\item  $X=\pa{X_{m}}_{m\ge 0}$,
$\widetilde{X}=\pa{\widetilde{X}_{m}}_{m\ge 0}$
   are two independent sequences with the same (unknown)
   distribution, independent of $(U,s_{\ell},W_{\vec\ell},s_{r},W_{\vec r})$,
  and therefore of $(A_{n,p},B_{n,p},Z_{n,p},T_{n,p})$.
\end{itemize}

The errors having a balancing effect:
$Z_{n,p}=\lceil nU \rceil-s_{\ell}+s_{r}$ has the
same mean, $(n+1)/2$, and a smaller variance than $\lceil nU \rceil$.
We prove this in the following form.

\begin{lem}\label{L:2}
\begin{align*}
\esp{(Z_{n,p}-1)^2 +(n-Z_{n,p})^2}
&\le
\esp{(\esup{nU}-1)^2+(n-\esup{nU})^2}
=\frac{(n-1)(2n-1)}{3} \\
&\le \frac23 n^2.
\end{align*}
\end{lem}
\begin{proof}
The left hand side is the expected number of ordered pairs $(i,j)$ that
end up on a common side of the pivot. This happens if $i$ and $j$
originally are on the same side of the pivot and we either compare
both correctly or make errors for both of them, or if they are on opposite
sides of the pivot and we make an error for exactly one of them.
Hence
\begin{align*}
\esp{(Z_{n,p}-1)^2 +(n-Z_{n,p})^2}
=\bigpar{p^2+(1-p)^2} \esp{(\esup{nU}-1)^2+(n-\esup{nU})^2}\quad&\\
  {}+2p(1-p)2 \esp{(\esup{nU}-1)(n-\esup{nU})}& \\
=\esp{(\esup{nU}-1)^2+(n-\esup{nU})^2}
-2p(1-p) \esp{(\esup{nU}-1-(n-\esup{nU}))^2}\qquad &
\end{align*}
which proves the first inequality. The rest is a simple calculation.
\end{proof}

Let us say that an element $a$ of the list,
 or the comparison in which $a$ plays the
 r\^ole of pivot,
 has depth $k$ if $a$ experiences 
$k-1$ comparisons before 
playing the r\^ole of pivot. We assume in this Section 
that  any comparison with depth $k+1$
is performed after the last comparison
with depth $k$.   We call  step $k$  the set of comparisons
with depth $k$, and 
we let $I^{(k)}(n,p)$ denote   
the number of inversions created at step $k$, that is, the total number 
of inversions, in the output,  between
  elements that are still in the same sublist before step $k$,
 but are not in the same sublist
after step $k$.
We shall need the following bound:
\begin{lem}\label{L:boundk}
For every $k\ge1$,
\[
\esp{I^{(k)}(n,p)}\le \frac12 \pa{\frac23}^k n^2p.
\]
\end{lem}

\begin{proof}
For $k=1$, $I^{(1)}(n,p)=t(n,p)$, and a simple calculation yields
\begin{align*}
\esp{t(n,p)}&= p\frac{(n-1)(n+1)}{3}-p^2\frac{(n-1)(n-2)}{6}
\le
  \frac13 n^2p.
\end{align*}

For $k>1$ we find by induction, conditioning on the partition in
the first step,
\[
\esp{I^{(k)}(n,p)}
\le
\esp{
\frac12 \pa{\frac23}^{k-1} (Z_{n,p}-1)^2p
+\frac12 \pa{\frac23}^{k-1}  (n-Z_{n,p})^2p
}
\]
and the result follows by Lemma \ref{L:2}.
\end{proof}

\begin{prop}
\label{bound}
Set $a_{n,p}=\esp{X_{n,p}}$. Then
\[a_{n,p}\le 1.\]
\end{prop}
\begin{proof}
By Lemma \ref{L:boundk},
$a_{n,p} \le \sum_1^\infty \frac12 \pa{\frac23}^k$.
\end{proof}

\section{Fixed point theorems}
\label{fixedpts}
The proofs of the first two cases are
  examples of the contraction method
\cite{Rosler91,RoeRue}:
  on one hand
we have more or less explicitly defined
random variables  $A^{(i)}_{n,p}$, $1\le i\le I$, and $T_{n,p}$, and we know
how to prove directly that they converge  to  $A^{(i)}$, $T$.
On the other hand, we have
a family $X_{n,p}$ of random variables  defined by induction:
\begin{align}
\label{masterequ1}
X_{n,p}&\overset{\text{law}}{=}\sum_{i=1}^I A^{(i)}_{n,p}
X^{(i)}_{Z^{(i)}_{n,p},p}+T_{n,p},
\end{align}
and a random variable $X$
 implicitly defined by
   the distributional identity
\begin{align}
\label{masterequ2}
X&\overset{\text{law}}{=}\sum_{i=1}^I A^{(i)}
X^{(i)}+T,
\end{align}
in which, in some sense,
   $\lim Z_{n,p}^{(i)}=+\infty$.  Then,
under additional technical conditions,
the  convergence of the "coefficients" $A^{(i)}_{n,p}$, $T_{n,p}$, entails
the convergence
of the "solution"  $X_{n,p}$. One has to prove  existence and unicity of
the solutions,
usually as fixed points of  contracting transformations in a subspace
of the space of
probability measures, with a suitable metric.
In the case we are interested in, (\ref{masterequ1}) holds and:
\begin{itemize}
\item
$I$ is a fixed positive integer;
\item
$C_{n,p}=(A_{n,p}^{(1)},Z_{n,p}^{(1)},
\dots,A_{n,p}^{(I)},Z_{n,p}^{(I)},T_{n,p})$
is a given random vector for each $n,p$;
\item
$Z_{n,p}^{(i)} \in [0,\dots,n-1]$;
\item
The families
$(X_{n,p}^{(i)})_{n,p}$, $i=1,2,\dots,I$, are
i.i.d.\ and independent of $C_{n,p}$, and
and $X_{n,p}^{(i)}\build{=}{}{law} X_{n,p}$.
\end{itemize}

Given such $C_{n,p}$ we thus define,
for any distributions
$G_{0,p},\dots,G_{n-1,p}$,
$$
\Phi(G_{0,p},\dots,G_{n-1,p}) =\lr\pa{\sum_{i=1}^I A^{(i)}_{n,p}
X^{(i)}_{Z^{(i)}_{n,p},p}+T_{n,p}},
$$
when, as above, the families
$(X_{k,p}^{(i)})_{k,p}$, $i=1,2,\dots,I$, are
i.i.d.\ and independent of $C_{n,p}$,
and further $X_{k,p}^{(i)}$ has the distribution $G_{k,p}$.
Thus (\ref{masterequ1}) can be written
\[G_{n,p}=\Phi(G_{0,p},\dots,G_{n-1,p}).\]

For (\ref{masterequ2}) we similarly assume
\begin{itemize}
\item
$C=(A^{(1)},\dots,A^{(I)},T)$
is a given random vector;
\item
the variables
$X^{(i)}$, $i=1,2,\dots,I$ are
i.i.d.\ and independent of $C$, and $X^{(i)} \build{=}{}{law} X$.
\end{itemize}
Given such $C$ we define
$$
\Psi(F)=\lr\pa{\sum_{i=1}^{I}A^{(i)}X^{(i)}+T},
$$
when
the variables
$X^{(i)}$, $i=1,2,\dots,I$ are
i.i.d.\ with distribution $F$ and independent of $C$.
Then (\ref{masterequ2}) can be written
\[\Psi(F)=F.\]

Let $D$ be the space of probability measures $\mu$ on
$\mathbb R$ such that $\int_{\mathbb R}\abs x\,d\mu(x)<+\infty$. The space
$D$ is endowed
   with the Wasserstein metric
\begin{align}
\label{metric}
d_1(\mu,\nu)&=
\inf_{{\scrip\lr(X)=\mu}\atop{\scrip\lr(Y)=\nu}}\left\|X-Y\right\|_1\\
\nonumber &=\left\|F^{-1}(U)-G^{-1}(U)\right\|_1.
\end{align}
in which $F$ and $G$ denote the distribution
functions of $\mu$ and $\nu$, $F^{-1}$ (resp.\ $G^{-1}$) denote the
  generalized inverses  of  $F$ and $G$  and, as in previous sections,
  $U$ is a uniform random variable \cite{Cambanis}. Since $F^{-1}(U)$
  (resp.\ $G^{-1}(U)$) has distribution $\mu$ (resp.\ $\nu$),
the infimum is attained  in relation (\ref{metric}).

The metric $d_1$ makes $D$  a complete metric space.
Convergence  of $\mathcal L(X_n)$ to  $\mathcal L(X)$ in $D$ is
equivalent to convergence of $X_n$ to $X$ in distribution \textit{and}
  \[\lim\esp{ |X_n|}=\esp{|X|}.\]
Therefore convergence   in $D$  entails
\[\lim\esp{X_n}=\esp{X}.\]
We refer to \cite{Rachev} for an extensive treatment of
Wasserstein metrics. In what  follows, we shall improperly
refer to the convergence  of $X_n$ to  $X$ in $D$,
meaning the convergence of their distributions.  Let us
take care first of relation (\ref{masterequ2}):

\begin{thm}\label{pointfixe}
If
  $\disp\sum_{i=1}^I\esp{\big|A^{(i)}\big|}<1$ and $\esp{|T|}<\infty$,
then
$\Psi$ is a strict contraction and (\ref{masterequ2}) has a unique
  solution in $D$.
\end{thm}

\begin{proof}
Let  $(X,Y)$ be a coupling
of  random variables, with laws $\mu$ and  $\nu$, respectively, such that
\[
\mathbb{E}\cro{\left\vert X-Y\right\vert}= d_1(\mu,\nu).
\]
Let  $\pa{\pa{X^{(i)},Y^{(i)}}}_{1\le i\le I}$
be $I$ independent copies of  $(X,Y)$. Furthermore,
assume that $C$ and $\pa{\pa{X^{(i)},Y^{(i)}}}_{1\le i\le I}$
are independent.  Then the probability distribution of
\[\sum_{i=1}^{I}A^{(i)}X^{(i)}+T,\quad  \textrm{resp.}\
\sum_{i=1}^{I}A^{(i)}Y^{(i)}+T\]
is $\Psi(\mu)$  (resp.\ $\Psi(\nu)$) and
\begin{eqnarray*}
d_1(\Psi(\mu),\Psi(\nu))
&\leq&\sum_{i=1}^I\esp{\big|A^{(i)}\big|\ \big|X^{(i)}-Y^{(i)}\big|}\\
&\le&d_1(\mu,\nu)\ \sum_{i=1}^I\esp{\big|A^{(i)}\big|}.
\end{eqnarray*}
Thus $\Psi$ is a contraction with contraction constant smaller than
1. Since $D$
  is a complete metric space,
  this implies that $\Psi$ has a unique  fixed point in $D$,
  by Banach's fixed point theorem.
\end{proof}

We prove now a theorem which is a variant
of those used by the previously cited authors: the difference is not deep,
but here we deal with family of laws, not sequences, as we have
  two parameters, $n$ and $p$.  As a consequence,
to cover Theorems \ref{p=c} and \ref{npinfini},
  it will be convenient  in their proofs to  consider convergence
with respect to a \textit{filter} $\mathcal F$ on $\mathbb N\times [0,1]$, see
\cite[Chap. 1]{Bourbaki}.
The collection of sets
\[
V_{N,\varepsilon }
=
\{n\ge N\}\times(\cro{c-\varepsilon,c+\varepsilon}\cap\left[0,1\right]),
\hspace{1cm}
N\ge 0,\ \varepsilon>0,
\]
is a basis for the filter $\mathcal F_1$  corresponding to Theorem \ref{p=c},
while 
\[\widetilde V_{N,\varepsilon }=\ac{(n,p)\,\big\vert\,0<p\le
\varepsilon,\,n\ge N/p},
\hspace{1cm}
N\ge 0,\ \varepsilon>0,
\]
is a basis for the filter $\mathcal F_2$  corresponding to Theorem \ref{npinfini}.

\begin{thm}
\label{th de convergence}
Suppose that  \eqref{masterequ1} holds for $n\ge1$ and $X_{0,p}=0$; i.e.\
$G_{n,p}=\Phi(G_{0,p},\dots,G_{n-1,p})$ for $n\ge1$ and
$G_{0,p}=\delta_0$, where $G_{n,p}=\mathcal L(X_{n,p})$.
If
\begin{itemize}
\item[i)]   $\pa{\esp{X_{n,p}}}_{n,p}$  is bounded,

\item[ii)] $\disp\sum_{i=1}^I\esp{\big|A^{(i)}\big|}<1,$

\item[iii)] $T_{n,p}\overset{L^1}{\underset{\mathcal F}{\To}}T,$\
$A^{(i)}_{n,p}\overset{L^1}{\underset{\mathcal F}{\To}}A^{(i)},$
\item[iv)] $\lim_{\mathcal F}\
\disp\esp{\big|A_{n,p}^{(i)}\big|\textrm{ ; }(Z_{n,p}^{(i)},p)\notin V}
=0,\qquad\forall V\in\mathcal F,$
\end{itemize}
then
$X_{n,p}$ converges in distribution to $F$, the 
 unique solution of the equation $\Psi(F)=F$
in $D$. More precisely,
$d_1(G_{n,p},F) \to 0$ along $\mathcal F$.
\end{thm}
We need a lemma before proving Theorem \ref{th de convergence}.
\begin{lem}\label{lemme rosler}
Assume that three families of nonnegative numbers
$(a_{n,p})_{0\le n,0<p<1}$, $(b_{n,p})_{0\le n,0<p<1}$,
  and
$\left(\gamma_{i,n,p},\,0\le n, 0\le i\le n, 0<p<1\right)$ satisfy the inequalities:
\[\ a_{n,p}\leq  b_{n,p}+\sum_{i=0}^{n-1}\gamma_{i,n,p}a_{i,p},
\]
Let $\mathcal F$ be a filter. Under the following assumptions:
\begin{itemize}
\item[--] $a_{n,p}$ is nonnegative and bounded,
\item[--]for some $\Gamma<1$ and some $V_0\in \mathcal F$,
$\forall (n,p)\in V_0,
\ \sum_{k=0}^{n-1}\gamma_{k,n,p}<\Gamma$,
\item[--]$\lim_{\mathcal F} b_{n,p}=0$,
\item[--]$\forall V\in {\mathcal F},
\ \lim_{\mathcal F}
\sum_{k\textrm{:}(k,p)\notin V} \gamma_{k,n,p}=0,$
\end{itemize}
we have
\[\lim_{\mathcal F}\ a_{n,p}=0.\]
\end{lem}

\begin{proof}[Proof of Lemma \ref{lemme rosler}.]
The proof is a variant of the proof  of \cite[Proposition 3.3]{Rosler91}.
Let  $M$ be a bound for $a_{n,p}$, and let
  \[ a=\limsup_{\mathcal F}\ a_{n,p}.\]
For any
$\epsilon>0$,   let $V_{\epsilon}\in \mathcal F$ be such that for
$(n,p)\in V_{\epsilon}$,
\[a_{n,p} \leq  a+\epsilon.\]
Then for $(n,p)\in V_{\epsilon}\cap V_0$ we have
\begin{eqnarray*}
a_{n,p} &\leq&\sum_{k\textrm{:}(k,p)\notin V_{\epsilon}}
\gamma_{k,n,p}a_{k,p}
+\sum_{k\textrm{:}(k,p)\in V_{\epsilon}}\gamma_{k,n,p}a_{k,p}+b_{n,p}\\
&\leq&M\ \sum_{k\textrm{:}(k,p)\notin V_{\epsilon}} \gamma_{k,n,p}+
(a+\epsilon)\Gamma+b_{n,p}.
\end{eqnarray*}
Taking lim sups, we obtain that for any $\epsilon>0$,
$$a\leq (a+\epsilon)\Gamma.$$
Thus $a\leq a\Gamma$, and so $a=0$.
\end{proof}

\begin{proof}[Proof of Theorem \ref{th de convergence}.]  We can
choose $X^{(i)}$  and
the family $\pa{X_{k,p}^{(i)}}_{0\le k,0<p<1}$  in such a way that
\[\esp{\Big|X_{k,p}^{(i)}
-X^{(i)}\Big|}=d_1(G_{k,p},F),\]
and we can also choose the families
$\pa{X^{(i)},(X_{k,p}^{(i)})_{k\ge 0}}_{0\le i\le I}$
to be i.i.d.
Then
\begin{align*}
d_1(G_{n,p},F)
&\leq\esp{\Big|\sum_{i=1}^{I}A_{n,p}^{(i)}X_{Z_{n,p}^{(i)},p}^{(i)}
+T_{n,p}-\sum_{i=1}^{I}A^{(i)}X^{(i)}-T\Big|}\\
&\leq\sum_{k=0}^{n-1}
\esp{\sum_{i=1}^{I}\big|A_{n,p}^{(i)}\BBone_{Z_{n,p}^{(i)}=k}\big|}
\esp{\big|X_{k,p}^{(i)}-X^{(i)}\big|}
+b_{n,p}\\
&\leq\sum_{k=0}^{n-1}\gamma_{k,n,p}\,d_1(G_{k,p},F)+b_{n,p}
\end{align*}
with
\begin{eqnarray*}
b_{n,p}&=&\sum_{i=1}^{I}\esp{\big|A^{(i)}_{n,p}-A^{(i)}
\big|X^{(i)}}+\esp{\abs{T_{n,p}-T}},\\
\gamma_{k,n,p}&=&\sum_{i=1}^{I}\esp{\big|A_{n,p}^{(i)}\big|
\BBone_{Z_{n,p}^{(i)}=k}}.
\end{eqnarray*}
Let  $M$ be a bound for $\pa{\esp{X_{n,p}}}_{n,p}$, and set
\[
a_{n,p}=
d_1(G_{n,p},F).
\]
Let us check the assumptions of Lemma \ref{lemme rosler}:
\[0\le a_{n,p}\le \esp{X_{n,p}}+\esp{X}\le M + \esp{X}\ ;\]
for the second assumption of Lemma \ref{lemme rosler},
$$\limsup_{\mathcal F}\sum_{k=0}^{n-1}\gamma_{k,n,p}
=\limsup_{\mathcal F}\sum_{i=1}^{I}\esp{\big|A^{(i)}_{n,p}\big|}
=\sum_{i=1}^{I}\esp{\big|A^{(i)}\big|}<1\ ;$$
  $\lim_{\mathcal F}b_{n,p}=0$ by assumption   iii), as
\begin{eqnarray*}
\sum_{i=1}^{I}\esp{\big|A^{(i)}_{n,p}-A^{(i)}
\big|X^{(i)}}&=&\esp{X}\sum_{i=1}^{I}\esp{\big|A^{(i)}_{n,p}-A^{(i)}\big|}\ ;
\end{eqnarray*}
finally
\[\sum_{k\textrm{ s.t. }(k,p)\notin V} \gamma_{k,n,p}=\sum_{i=1}^{I}
\disp\esp{\big|A_{n,p}^{(i)}\big|\textrm{ ; }(Z_{n,p}^{(i)},p)\notin V}.\]
Therefore $d_1(G_{n,p},F)$ vanishes along $\mathcal F$ and the proof of
the theorem is now complete.
\end{proof}

The following Theorem is folklore. It gives  the means and variances
in Theorems \ref{p=c} and \ref{npinfini},
after some computations.
\begin{thm}
\label{usefulthmaboutvariance}
Suppose that \eqref{masterequ2} holds, where
$\sum_i \esp{|A^{(i)}|}<1$  and $\esp{|X|}<\infty$;
in other words,
$\mathcal L(X)=F$, where
$F$ is the unique solution in $D$
to $\Psi(F)=F$.
Then
\begin{equation}
\label{noel}
\esp{X} = \frac{\esp{T}}{1-\sum_i \esp{A^{(i)}}}.
\end{equation}

Moreover, if further
$\sum_i \esp{|A^{(i)}|^2}<1$ and $\esp{ T^2}<\infty$,
then  $\esp{X^2}<\infty$ and
\begin{equation}
\label{jul}
\var{X} = \frac{\esp{T^2}+ 2\,\esp X \esp{T\sum_i A^{(i)}}
+\esp{\pa{\sum_i A^{(i)}}^2-1}(\esp X)^2}
{1-\sum_i \esp{A^{(i)\,2}}}.
\end{equation}
\end{thm}

\begin{proof}
Taking expectations in (\ref{masterequ2}) we obtain
$\esp X= \sum_i \esp{A^{(i)}} \esp X + \esp T$,
which yields (\ref{noel}).

For the second part, let
$D_2=\{\mu\in D: \int x\,d\mu(x)=\esp X,\, \int x^2\,d\mu(x)<\infty\}$.
It is easy to see that now $\Psi$ is a strict contraction in $D_2$
with the $d_2$ metric; hence $\Psi$ has a unique fixed point in
$D_2$. Since $D_2\subset D$, this fixed point must be $F$, which shows
that
$\esp{X^2} <\infty$.
If we square \eqref{masterequ2} and take the expectation, we obtain
\begin{align*}
\esp{X^2}
&= \esp{\sum_{i=1}^I A^{(i)\,2}} \esp{X^2}
+\sum_{1\le i\neq j\le I} \esp{A^{(i)}A^{(j)}} (\esp X)^2\\
&\hspace{3.5cm}+ 2 \sum_{i=1}^I \esp{A^{(i)}T} \esp X
+ \esp{ T^2},
\end{align*}
which yields \eqref{jul}.
\end{proof}

\section{Proofs of  Theorems \ref{p=c} and \ref{npinfini}}
\label{proofTh1}
We apply Theorem \ref{th de convergence} to the distributional identity
(\ref{eq en X}), with $I=2$,
\[
\pa{A_{n,p}^{(1)},Z_{n,p}^{(1)}}=\pa{A_{n,p},Z_{n,p}-1},
\]
and
\[
\pa{A_{n,p}^{(2)},Z_{n,p}^{(2)}}
=
\pa{B_{n,p},n-Z_{n,p}}.
\]
Here the distribution of $\pa{A_{n,p}^{(i)},Z_{n,p}^{(i)}}$ does not
depend on $i$.
We verify the assumptions ii)--iv) of
Theorem \ref{th de convergence} for
Theorems \ref{p=c} and \ref{npinfini} together;
for the second theorem take $c=0$.
The first assumption holds true  by Proposition \ref{bound}.

\subsection{Verification of the second point}
We have
\begin{align*}
A^{(1)}&=A=[(1-2c)U+c]^2,\\
A^{(2)}&=B=[(2c-1)U+1-c]^2,
\end{align*}
and $c\in[0,1]$.
Easy computations give
$$\esp{[(1-2c)U+c]^2}+\esp{[(2c-1)U+1-c]^2}=\frac23(1-c+c^2)\le \frac 23.$$

\subsection{Verification of the third point}

We must prove the convergence of $A_{n,p}$, $B_{n,p}$ and $T_{n,p}$ to
$A$, $B$ and $T(c,U)$,  in
$L^1$.
Recall \eqref{znp}--\eqref{tnp}.

From Proposition \ref{description}
 we know that, conditioned on $U$, $s_\ell\sim\Bi(\esup{nU}-1, p)$
and thus
\begin{equation*}
\mathbb{E}\bigpar{(s_\ell-(\esup{nU}-1)p)^2\mid U}
=(\esup{nU}-1)p(1-p) \le np.
\end{equation*}
Hence, taking the expectation,
\begin{equation*}
\mathbb{E}\bigpar{s_\ell-(\esup{nU}-1)p}^2
\le np
\end{equation*}
and thus
\begin{equation*}
\norm{s_\ell-nUp}_2
\le \norm{s_\ell-(\esup{nU}-1)p}_2 +p
\le (np)^{1/2} + p
\le 2(np)^{1/2}.
\end{equation*}
Consequently,
\begin{equation}\label{sl}
\norm{\frac{s_\ell}{n}-Uc}_2
\le
\norm{\frac{s_\ell}{n}-Up}_2+|p-c|
\to0,
\end{equation}
and, similarly but more sharply,
\begin{equation}\label{slp}
\norm{\frac{s_\ell}{n\sqrt p}-U\sqrt c}_2
\le \frac 2{\sqrt n}+|\sqrt p-\sqrt c|
\to0.
\end{equation}
   Similarly,
\begin{equation}\label{sr}
\norm{\frac{s_r}{n}-(1-U)c}_2
\to0
\end{equation}
and
\begin{equation}\label{srp}
\norm{\frac{s_r}{n\sqrt p}-(1-U)\sqrt c}_2
\to0.
\end{equation}

{}From \eqref{znp}, \eqref{sl} and \eqref{sr} follows
\begin{equation}\label{sofie}
\norm{\frac{Z_{n,p}-1}{n}-\bigpar{U-Uc+(1-U)c}}_2
\to0.
\end{equation}
It follows easily from Cauchy--Schwarz's inequality that multiplication is a
continuous bilinear map $L^2\times L^2\to L^1$. Hence
\eqref{sofie} yields
\begin{equation*}
\norm{A_{n,p}-A}_1=
\norm{\pa{\frac{Z_{n,p}-1}{n}}^2-\bigpar{U-Uc+(1-U)c}^2}_1
\to0,
\end{equation*}
verifying the first assertion.
\eqref{sofie} similarly implies
$\norm{B_{n,p}-B}_1\to0$ too.

For $T_{n,p}$ we first observe that, similarly, from \eqref{slp} and
\eqref{srp},
\begin{equation*}
\norm{\frac{s_\ell s_r}{n^2p}-U(1-U)c}_1
\to0.
\end{equation*}
Moreover, since $np\to\infty$, \eqref{sl} and \eqref{sr} imply
$\norm{\xfrac{s_\ell}{n^2p}}_1
\le \norm{\xfrac{s_\ell}{n^2p}}_2 \to0$
and $\norm{\xfrac{s_r}{n^2p}}_1\to0$.

For the terms $W_{\vec\ell}$ and $W_{\vec r}$ we use Proposition
\ref{invb&w}.
We have
$\norm{S_m-mp}_2=\sqrt{mp(1-p)}$
and thus, uniformly for $0\le m\le n$,
\begin{equation*}
\norm{\frac{S_m}{n\sqrt p} -\frac{m}{n}\sqrt c}_2
\le
\frac{1}{\sqrt n}+ |\sqrt p - \sqrt c|
\to 0,
\end{equation*}
which, using Cauchy--Schwarz again, yields
\begin{equation}\label{erika}
\norm{\frac{S_m(S_m+1)}{2n^2 p} -\frac{c}{2}\pa{\frac{m}{n}}^2}_1
\to 0.
\end{equation}
Moreover, let $W'_m=\sum_{i=1}^m i Y_i$. Then
$\mathbb{E} W'_m = \frac{m(m+1)p}2$ and
\begin{equation*}
\norm{W'_m-\mathbb{E} W'_m}_2^2 =\var{W'_m} = \sum_{i=1}^m i^2p(1-p) \le
m^3p,
\end{equation*}
and thus
\begin{equation}\label{magnus}
\norm{\frac{W'_m}{n^2p} -\frac{1}{2}\pa{\frac{m}{n}}^2}_2
\le\frac{1}{\sqrt{np}}+\frac{1}{2n} \to 0.
\end{equation}

Proposition \ref{invb&w} now yields, by \eqref{erika} and
\eqref{magnus},
uniformly for $m\le n$,
\begin{equation*}
\norm{\frac{W_m}{n^2p} -\frac{1-c}{2}\pa{\frac{m}{n}}^2}_1 \to 0.
\end{equation*}
Consequently, using
Proposition \ref{description},
\begin{align*}
\norm{\frac{W_{\vec\ell}}{n^2p} -\frac{1-c}{2}U^2}_1 &\to 0, \\
\norm{\frac{W_{\vec r}}{n^2p} -\frac{1-c}{2}(1-U)^2}_1& \to 0.
\end{align*}
Collecting the various terms above,
we find $\norm{T_{n,p}-T}_1\to0$.

\subsection{Verification of the fourth point}

As already noticed at the beginning of the Section, 
the distribution of $\pa{A_{n,p}^{(j)},Z_{n,p}^{(j)}}$ does not
depend on $j\in\{1,2\}$, so in order to prove the two theorems, we only 
 have  to check that the fourth assumption holds for $j=1$, 
for an arbitrary set in each of the  two filters: 
\begin{equation}
\label{equation_filter}
\lim_{\mathcal F_i}\ \esp{\big|A_{n,p}\big|\textrm{ ; }(Z_{n,p}-1,p)\notin V}
=
0,
\qquad
\forall V\in\mathcal F_i, \forall i\in\{1,2\} ;
\end{equation}
also,  the  expectation on the left hand side of (\ref{equation_filter})
 is decreasing in $V$, so 
we need only to check (\ref{equation_filter}) for  typical elements
of the filters' basis.
But  for $(n,p)\in V_{N,\varepsilon}$  (resp.\  for $(n,p)\in
  \widetilde V_{N,\varepsilon}$),
\begin{align*}
\esp{\big|A_{n,p}\big|\textrm{ ; }(Z_{n,p}-1,p)\notin V_{N,\varepsilon}}
&\le \pa{\frac{N-1}n}^2,\\
\esp{\big|A_{n,p}\big|\textrm{ ; }(Z_{n,p},p)\notin
  \widetilde V_{N,\varepsilon}}
&\le \pa{\frac{N}{np}}^2.
\end{align*}

\section{Proofs of Theorems \ref{p=1/n} and \ref{thm:unique}}
\label{proofTh3}

The proof of these theorems is done in four steps:
\begin{enumerate}
\item
We prove   that  $X(\gl)$  defined at (\ref{def de I})
is almost surely finite, and has
exponentially decreasing distribution tail.
Thus it has moments of all orders.
\item
With the help
of a Poisson point process representation of Quicksort,
  we prove the convergence of certain
copies of $X_{n,p}$ to a copy of $X(\gl)$
  for the norm $\norm{\cdot}_1$.  This entails the weak convergence.
\item
We prove that  $X(\gl)$ satisfies the functional
   equation   (\ref{functional1}), and that (\ref{functional1}) has a
unique solution under the extra assumptions in Theorem \ref{thm:unique}.
\item
We compute the first and second moments of  $X(\gl)$, as
 required for  the proof of  Theorem \ref{p=1/n}, and 
we also give an induction formula for moments of larger order.
\end{enumerate}

\subsection{Some properties of $X(\gl)$.}
In this Section, we prove some properties of the family of random variables
$\pa{X(\gl)}_{\gl>0}$ defined by   \eqref{def de I}.
Recall that the increasing sequence $\pa{Y_{k,j}}_{0\le j\le 2^k}$, defined
  by the  recurrence relation (\ref{recurrence_sur_U}),
splits $[0,1]$ in  $2^k$ intervals, obtained  recursively
  by breaking each of the $2^{k-1}$ intervals of the previous step
  into two random pieces. For $k\ge 0$ and $1\leq  i\leq 2^k$, let
\begin{align*}
w_{k,i}&=Y_{k,i}-Y_{k,i-1},\\
M_k&=\max\Big\{w_{k,i}: 1\leq  i\leq 2^k\Big\},\\
F_{k,\alpha}&=\pa{\frac{1+\alpha}2}^k\sum_{1\leq  i\leq 2^k}w_{k,i}^{\alpha}\\
\mathcal{F}_k&=\sigma\pa{Y_{i,j}:i\le k,1\leq j\leq2^i-1}\\
\mathcal{F}&=\pa{\mathcal{F}_k}_{k\ge 0}.
\end{align*}
We begin with a simple estimate (see also \cite{Gruber_Rosler96}):
\begin{lem}
\label{L:erika}
$\esp{w_{k,j}^2}= 3^{-k}$.
\end{lem}
\begin{proof}
The length $w_{k,j}=|\ikj|$ is the product of $k$ independent random
variables, each
uniform on $[0,1]$. Hence $\esp{w_{k,j}^2}=\pa{\esp{U^2}}^k=3^{-k}$.
\end{proof}

\begin{lem}\label{Bigg} For $\alpha > 0$,
$\pa{F_{k,\alpha}}_{k\ge 0}$
is a $\mathcal{F}$-martingale, and $\esp{F_{k,\alpha}}=1$.
\end{lem}

\begin{proof} Clearly $\esp{F_{0,\alpha}}=1$. Also:
\begin{align*}
\esp{F_{k+1,\alpha}|\mathcal{F}_k}
&=\pa{\frac{1+\alpha}2}^{k+1}\sum_{i=1}^{2^k}\esp{w^{\alpha}_{k+1,2i-1}
+w^{\alpha}_{k+1,2i}|\mathcal{F}_k}\\
&=\pa{\frac{1+\alpha}2}^{k+1}\sum_{i=1}^{2^k}w^{\alpha}_{k,i}
\esp{{U}_{k,i}^{\alpha}+(1-{U}_{k,i})^{\alpha}}\\
&=\pa{\frac{1+\alpha}2}^{k}\sum_{i=1}^{2^k}w^{\alpha}_{k,i}.
\end{align*}
\end{proof}
Let $\rho=0.792977\dots$
denote the larger real solution of the equation
$\rho^{-1}=-2e\ln \rho$.  Lemma \ref{Bigg} entails  that
\begin{lem}
\label{L:moment}
$\esp{M_k}\le \rho^{k}$.
\end{lem}
\begin{proof}Clearly,
\[M_k^{\alpha}\leq \pa{\frac 2{1+\alpha}}^k F_{k,\alpha};\]
  thus, for $\alpha\ge 1$,
\[\esp{M_k}\le\left(\esp{M_k^{\alpha}}\right)^{1/\alpha}\leq
\pa{\frac2{1+\alpha}}^{k/\alpha}\left(\esp{F_{k,\alpha}}\right)^{1/\alpha}
=\pa{\frac2{1+\alpha}}^{k/\alpha}.
\]
The rate $\pa{\frac2{1+\alpha}}^{1/\alpha}$ reaches its minimum
  for $1+\alpha=4.311\dots$, a constant that is an old friend  of Quicksort
and binary search trees \cite{Devroye}. This leads to
the desired value for $\rho$.
\end{proof}

A  weaker form of this inequality (for $\alpha=2$), actually sufficient
for our purposes, is given in \cite{Gruber_Rosler96}.
The sequence $\pa{F_{k,\alpha}}_{k\ge 0}$
  is a specialization  of martingales that are of a great use for the study of
general branching random walks, see for instance \cite{Biggins2},
of which binary search trees are  a special case \cite{Mahmoud,Pittel}.

\begin{lem}
\label{L:mean}
$\esp{X(\gl)}=1$.
\end{lem}
\begin{proof}
Set
$\mathcal{F}_{\infty}=\sigma\pa{Y_{k,j},k\ge 0,1\leq j\leq2^k-1}$.
Inspecting \eqref{emma}, we see that
\[\esp{X(\gl)\,\vert\,\mathcal{F}_{\infty}}=\frac1 2\  \sum_{k\ge 1}\pa{\frac23}^{k}F_{k,2},\]
because, conditionally given $\mathcal{F}_{\infty}$, the expected number of points of $\Pi_{k,j}$
is $\gl w_{k,j}$ and each of them has an expected contribution $w_{k,j}/(2\gl)$ to
$X(\gl)$.
\end{proof}
As a     consequence  of  Lemma \ref{L:moment}, we have
\begin{thm}
\label{tail}
For each fixed $\lambda>0$, the distribution tail $\pr{X(\gl)\geq x}$ decreases exponentially fast.
\end{thm}
\begin{proof}
Equivalently, we prove this result for $\Xi(\lambda)=\gl X(\gl)$. Since
\[|x-Y_{k,J_k(x)}\big|
\le
M_k,\]
we have
\begin{align*}
\Xi(\lambda)&\le\sum_{(k,x)\in\Pi} M_k =\sum_{k\ge 1} N_kM_k,
\end{align*}
where
$
N_k=\left\vert\Pi_k\right\vert
$
is a Poisson random variable with mean $\gl$. We split the tail of this bound on 
$\Xi(\lambda)$ as follows:
\begin{align*}
\pr{\Xi(\lambda)\ge x}&\le\pr{\sum_{k\ge 1} N_kM_k\ge x}
\\
&\le p_1+p_2,
\end{align*}
in which
\begin{align*}
p_1&=\pr{\sum_{1\le k\le m} N_kM_k\ge x/2},
\\
p_2&=\pr{\sum_{k> m} N_kM_k\ge x/2}.
\end{align*}
We have, by the standard Chernoff bound for the Poisson distribution,
\[
p_1\le\pr{\sum_{1\le k\le m} N_k\ge x/2}
\le\exp\pa{\frac x2\pa{1-\ln(x/2m\lambda)}-n\gl},
\]
the last inequality holding only for $m\le \frac x{2\gl}$. Also
\begin{align*}
p_2 &\le
\pr{\sum_{k> m} N_k\rho^{k/2}\ge x/2}
+\pr{\exists k>m: \, M_k>  \rho^{k/2}}
\\
&\le
\pa{1+\frac{2\lambda}x}\frac{\rho^{(m+1)/2}}{1-\sqrt \rho},
\end{align*}
using a Markov first moment inequality to bound both terms. For any
$\alpha$ in $(0,1)$,
the choice  $m\sim \frac{\alpha x}{2\gl}$ leads to an
exponential decrease of the tail.
\end{proof}

\subsection{Convergence of $X_{n,p}$ to  $X(\gl)$.}
We assume that the input list for Quicksort
contains the integers $\{1,2,\dots,n\}$
in random order. We model our error-prone Quicksort as follows using the
variables
$U_{k,j}$ and $\Pi$ in Section  \ref{S:results},
but with the intensity
$\gl$ of $\Pi$ replaced by $\gl(n,p)=-n\ln(1-p)$:

In the first step, we use the pivot
$p_{1,1}=\esup{n U_{0,1}}$ and let for each $i$ (except the pivot)
there be an error in the comparison of $i$ and the pivot if
$\Pi_1\cap(\frac{i-1}n,\frac{i}n] \neq\emptyset$.
(Recall that $\Pi_k=\set{x:(k,x)\in\Pi}$.)
Note that our choice of $\gl(n,p)$ yields the right error probability $p$.

Let $p'_{1,1}$ be the position of the pivot after the first step.
(This position was earlier denoted $Z_{n,p}$; it may differ from
$p_{1,1}$ because of errors.)
The items of the left sublist will thus be placed in positions
  $1,\dots,p'_{1,1}-1$ and those in
the right sublist in positions $p'_{1,1}+1,\dots,n$.
Let $p'_{1,0}=0$ and $p'_{1,2}=1+n$.

When the $k$-th step begins, we have a set of $2^{k-1}$ sublists
$\pa{\ell_{k-1,j}}_{j=1,\dots,2^{k-1}}$, the elements of $\ell_{k-1,j}$
being  in positions
$p'_{k-1,j-1}+1,\dots,p'_{k-1,j}-1$, $j=1,\dots,2^{k-1}$
(with the convention that the sublist is empty
when $p'_{k-1,j} - p'_{k-1,j-1} \le 1$).
In each nonempty such sublist we choose as pivot the item with rank
  $
\esup{U_{k-1,j}(p'_{k-1,j}-p'_{k-1,j-1}-1)},
$
in this sublist, so that its position in the final output will be exactly
\begin{equation}
\label{fauxpivot}
p_{k,2j-1}=p'_{k-1,j-1} + \esup{U_{k-1,j}(p'_{k-1,j}-p'_{k-1,j-1}-1)},
\end{equation}
in case  no errors occurs while processing the sublist.
We assume an error is made when comparing the element at position $i$ with the
pivot  $p_{k,2j-1}$ if
$\Pi_k\cap(\frac{i-1}n,\frac{i}n] \neq\emptyset$.
Let $p'_{k,2j}=p'_{k-1,j}$.  Let $p'_{k,2j-1}$ be the position of
  the pivot $p_{k,2j-1}$ after the comparisons (as in the first step,
$p'_{k,2j-1}$ may differ from  $p_{k,2j-1}$ because of errors);
let $p'_{k,2j-1}=p'_{k-1,j}$ if the sublist was empty.
Set
  \[y_{k,j}=p_{k,j}/n\hspace{0.3cm}\text{and}\hspace{0.3cm}y'_{k,j}=p'_
{k,j}/n.\]
   We expect $y_{k,j}$  and $y'_{k,j}$
to converge to $Y_{k,j}$ as $n\to +\infty$.

This procedure (stopped when there are no more nonempty sublists)
is an exact simulation of the erratic Quicksort, so we may assume
that $I(n,p)$ is the
number of inversions created by it.
As in Section \ref{functequns},
let $I^{(k)}(n,p)$ be the number of inversions
created at step $k$, so
\[I(n,p)=\sum_{k=1}^\infty I^{(k)}(n,p).\]

We will prove that, using the notation of \eqref{emma},
\begin{equation}
\label{jesper}
\gd_{n,k}= \norm{\frac 1{n^2p} I^{(k)}(n,p)-\frac1{\gl(n,p)}
\sum_{j=1}^{2^k} \sum_{x\in\Pi_{k,j}}|x-x_{k,j}|}_1
\To 0
\end{equation}
for each $k$.
Since also, by Lemmas \ref{L:boundk} and \ref{L:erika},
\begin{align*}
\gd_{n,k}
&\le
\frac 1{n^2p} \esp{I^{(k)}(n,p) }
+ \frac1{\gl(n,p)} \esp{\sum_{j=1}^{2^k} \sum_{x\in\Pi_{k,j}}|x-x_{k,j}|}\\
&=\frac 1{n^2p} \esp{I^{(k)}(n,p) }
+ \esp{ \frac12  \sum_{j=1}^{2^k}w_{k,j}^2}
\le \pa{\frac23}^k,
\end{align*}
it follows by dominated convergence that, using \eqref{emma},
\begin{equation*}
\norm{X_{n,p}-X\bigpar{\gl(n,p)}}_1
\le \sum_{k=1}^\infty \gd_{n,k} \To 0.
\end{equation*}
Moreover, $\gl(n,p)\to\gl$, and it follows easily from \eqref{samuel}
that
$
\norm{X\bigpar{\gl(n,p)}-X(\gl)}_1\to0
$.
Hence we have
$
\E\abs{X_{n,p}-X(\gl)} \to0
$, which proves the convergence.

It remains to verify \eqref{jesper}. Set
\[X^{(k)}=\sum_{j=1}^{2^k} \sum_{x\in\Pi_{k,j}}|x-x_{k,j}|.\]
Relation \eqref{jesper} is  equivalent to
\begin{equation}
\label{jesper2}
\E \abs{ \frac1n I^{(k)}(n,p)-X^{(k)}}\to 0.
\end{equation}
For simplicity, we write in the sequel $\gl$ instead of $\gl(n,p)$.
We begin with a lemma.

\begin{lem} \label{LL}
For each $k$ and $j$,
\[
  \max\ac{ \norm{Y_{k,j} - y'_{k,j}}_1,
\norm{Y_{k,j} -y_{k,j}}_1} \le \frac{k(1+\gl)}{n}.
\]

\end{lem}
\begin{proof}
Recall that $p'_{k,j}=ny'_{k,j}$, so
\eqref{fauxpivot} translates to
\begin{equation*}
p_{k,2j-1}=ny'_{k-1,j-1} + \esup{U_{k-1,j}(ny'_{k-1,j}-ny'_{k-1,j-1}-1)},
\end{equation*}
   We use induction on $k$.
Comparing the definitions of $Y_{k,j}$ and $y'_{k,j}$, we see that it
suffices to consider an  odd $j=2l-1$, and in that case there are three sources
of a difference:
\begin{enumerate}
\item The differences between $y'_{k-1,l-1}$ and $Y_{k-1,l-1}$
and between $y'_{k-1,l}$ and $Y_{k-1,l}$. By the induction hypothesis,
this contributes at most $(k-1)(1+\gl)/n$.
\item
The $-1$ inside  (and the rounding by) the ceiling function. This
contributes at most $1/n$.
\item
The shift of the pivot, from $p_{k,2j-1}$ to $p'_{k,2j-1}$,
caused by the erroneous comparisons.
The shift is bounded by the total number of errors at step $k$, so its
mean is less than $\gl$, and the contribution is less than $\gl/n$.
\end{enumerate}
\end{proof}

We return to proving \eqref{jesper2}.
For $k=1$, $I^{(1)}(n,p)$ is just $t(n,p)$ studied in Section
\ref{functequns}, and
\eqref{t1np}
yields
$$
I^{(1)}(n,p)=s_{\ell}s_{r}+s_{\ell}+s_{r}+W_{\vec\ell}+W_{\vec r} .
$$
Let $E_1$ be the set of items $i$ such that an error was made in the
comparison with $p_{1,1}$.
Relation \eqref{poissonstep1}  entails
that
$$
\sum_{i\in E_1}|i-p_{1,1}|
= W_{\vec\ell}+W_{\vec r}
+\tfrac12s_\ell(s_\ell+1)+\tfrac12 s_r(s_r+1).
$$
We shall denote this last sum $\tilde I^{(1)}(n,p)$. Thus, we have
$$
\abs{I^{(1)}(n,p)-\tilde I^{(1)}(n,p)}
=\abs{s_{\ell}s_{r}+s_{\ell}+s_{r}
-\tfrac12 s_\ell(s_\ell+1)-\tfrac12 s_r(s_r+1)}
\le s_\ell^2+s_r^2.
$$
Furthermore
\begin{align}
\esp{s_{\ell}^2\mid\esup{nU_{1,1}}=m}&=(m-1)p(1-p)+((m-1)p)^2\le np+n^2p^2.
\label{boundagain}
\end{align}
Hence,
$$
\norm{I^{(1)}(n,p)-\tilde I^{(1)}(n,p)}_1
= \bigO{1}.
$$
Moreover, $\frac1n\tilde I^{(1)}(n,p)=\sum_{i\in E_1} |\frac in-y_{1,1}|$
differs from
$X^{(1)}=\sum_{j=1}^{2} \sum_{x\in\Pi_{1,j}}|x-x_{1,j}|$ in \eqref{jesper2}
in four ways only
(recall that $x_{1,1}=x_{1,2}=Y_{1,1}$):
\begin{enumerate}
   \item
$i/n$ differs from $x$ by at most $1/n$. Since the expected number of
   terms is not larger than $\gl$, this gives a contribution $\bigO{1/n}$.
\item
$|y_{1,1}-x_{1,j}|=|y_{1,1}-Y_{1,1}|$,
which by Lemma \ref{LL} has expectation $\bigO{1/n}$.
Thus this too gives a contribution $\bigO{1/n}$.
\item
If there are two or more points in $\Pi_1\cap(\frac{i-1}n,\frac{i}n]$
for some $i$, $X^{(1)}$  contains more terms than
$\frac1n\tilde I^{(1)}(n,p)$.
It is easily seen that the expected number of such extra points in each
interval $(\frac{i-1}n,\frac{i}n]$  is less  than $(\gl/n)^2$,
and each point contributes at most 1 to $X^{(1)}$.
\item
Each point in $\Pi_1\cap(\frac{p_{1,1}-1}n,\frac{p_{1,1}}n]$
contributes for an extra term in $X^{(1)}$  again.
The expected number of such extra points is $\gl/n$
and each of these terms contributes at most 1 to $X^{(1)}$.
\end{enumerate}
This verifies \eqref{jesper2} for $k=1$.

For $k\ge2$ we argue similarly.
We can approximate $I^{(k)}(n,p)$ by the sum
of the distances between
the errors and the respective pivots,
\begin{align*}
\tilde I^{(k)}(n,p)
&=\sum_{j\le 2^{k-1}}\sum_{i\in  E_{k-1,j}}|i-p_{k,2j-1}|,
\end{align*}
as follows:  Let $E_{k,j}$ be the set of
items $i\in\ell_{k,j}$  subject to error when
compared with $p_{k+1,2j-1}$, and let $\mathcal G_k$ be the $\sigma$-algebra
generated by $\pa{U_{\ell,j}}_{\ell\le k,j\le 2^\ell} $ and
$\Pi_1\cup \Pi_2\cup \dots\cup \Pi_{k-1}$.
As for $k=1$, using relation \eqref{poissonstep1}, we obtain
the following bound:
\begin{align*}
\esp{\left.\left\vert I^{(k)}(n,p)-\tilde I^{(k)}(n,p)
\right\vert\ \ \right\vert\ \ \mathcal G_k}
&\le \sum_{j\le 2^{k-1}}\pa{p^2\pa{\#\ell_{k-1,j}}^2+p\#\ell_{k-1,j}}\\
&\le 2^{k-1}(n^2p^2+np)
=\bigO{1},
\end{align*}
and as a consequence,
$$
\norm{I^{(k)}(n,p)-\tilde I^{(k)}(n,p)}_1
=\bigO{1}.$$

Now,
\[\frac1n\tilde I^{(k)}(n,p)=\sum_{j\le 2^{k-1}}\sum_{i\in  E_{k-1,j}}
\left\vert\frac in-y_{k,2j-1}\right\vert\]
differs from
$X^{(k)}=\sum_{j=1}^{2^k} \sum_{x\in\Pi_{k,j}}|x-x_{k,j}|$ in \eqref{jesper2}
in the same four ways as for $k=1$, plus an extra fifth way:
\begin{enumerate}
   \item
See the case  $k=1$.
\item
$|y_{k,2j-1}-x_{k,2j-1}|=|y_{k,2j-1}-x_{k,2j}|=|y_{k,2j-1}-Y_{k,2j-1}|$,
which by Lemma \ref{LL} has expectation $\bigO{1/n}$.
Thus this too gives a contribution $\bigO{1/n}$.
\item
Two or more points in $\Pi_k\cap(\frac{i-1}n,\frac{i}n]$
for some $i$, see the case  $k=1$.
\item
Each point in $\Pi_k\cap(-1/n+y_{k,2j-1},y_{k,2j-1}]$
contributes for an extra term in $X^{(k)}$.
The expected number of such extra points is $\gl2^{k-1}/n$
and each of these terms contributes at most 1 to $X^{(k)}$.
\item
There is a new source of error in this approximation, because
some points $x$ in $\Pi_k$ and the corresponding positions $i=\esup{nx}$
  belong to subintervals that do not correspond to each
other, because the endpoints $y'_{k-1,j}$ differ somewhat from
$Y_{k-1,j}$.
By Lemma \ref{LL}, the expected number of such cases is $\bigO{1/n}$, so
again we get a contribution of order $\bigO{1/n}$ only.
\end{enumerate}

This verifies \eqref{jesper2} and thus the
convergence of $X_{n,p}$ to  $X(\gl)$.

\subsection{The distributional identity  for $X(\gl)$.}

We  check that $X(\gl)$ satisfies the distributional identity and 
some side conditions needed for the computations of moments.
\begin{prop}
\label{ilexiste1}
$\pa{X(\gl)}_{\gl>0}$  is a
solution of \eqref{functional1}.
Moreover, $\esp{X(\gl)^n}<\infty$ and
$\gl^n\esp{X(\gl)^n}\to0$ as $\gl\to0$, for $n\ge1$.
\end{prop}
\begin{proof}
All moments are finite by Theorem \ref{tail}.
Moreover, $\esp{(\gl X(\gl))^n}\to0$ as $\gl\to0$ by \eqref{samuel}
and dominated convergence.

For $a<b$, let $\Pi(a,b)$ be a Poisson point process of intensity $\gl$ on
${\mathbb N}^*\times[a,b]$,
and let $\{U_{k,j}: k\ge 0, 1\le i\le 2^k\}$ be
independent uniform random variables
as in Section \ref{S:results},
and further independent of $\Pi(a,b)$.
Define
$\{Y_{k,j}: k\ge 0, 1\le i\le 2^k\}$
and $J_k(x)$
as in Section \ref{S:results}, with the slight modification
\[Y_{0,0}=a\qquad\text{and}\qquad Y_{0,1}=b\]
  and set
$$
X(\gl,a,b)= \frac1{\gl}\ \sum_{(k,x)\in\Pi(a,b)}\big|x-Y_{k,J_k(x)}\big|.
$$
Note that $X(\gl,0,1) = X(\gl)$. Shifting and rescaling
$\Pi(a,b)$, we obtain
\[X(\gl, a, b)\egenloi X(\gl, 0, b-a)\egenloi(b-a)^2X(\gl(b-a)).\]
Let us split $X(\gl)$: we have
\begin{align*}
X(\gl)&=X_0(\gl)+X_1(\gl)+X_2(\gl)
\\
\gl X_0(\gl)&=\sum_{(1,x)\in\Pi(0,1)}\big|x-Y_{1,1}\big|
\\
\gl X_1(\gl)&=\sum_{(k,x)\in\Pi(0,1)\atop{k\ge 2,x\le Y_{1,1}}}
\big|x-Y_{k,J_k(x)}\big|,
\\
\gl X_2(\gl)&=\sum_{(k,x)\in\Pi(0,1)\atop{k\ge 2, x\ge Y_{1,1}}}
\big|x-Y_{k,J_k(x)}\big|.
\end{align*}
We see, using general properties of Poisson point processes  and the recursive
construction of $\{Y_{k,j}: k\ge 0, 1\le i\le 2^k\}$, that
\begin{align*}
(X_0(\gl),X_1(\gl),X_2(\gl))&\egenloi
(\Theta(\gl,Y_{1,1}),X(\gl,0,Y_{1,1}),\widetilde X(\gl,Y_{1,1},1))\\
&\egenloi(\Theta(\gl,Y_{1,1}),
Y_{1,1}^2X(\gl Y_{1,1}),(1-Y_{1,1})^2\widetilde X(\gl(1-Y_{1,1}))),
\end{align*}
in the sense
that, conditionally given that $Y_{1,1}=u$,
$X_0(\gl)$, $X_1(\gl)$ and $X_2(\gl)$
are independent and distributed as
$\Theta(\gl,u)$, $u^2X(\gl u)$, $(1-u)^2 X(\gl(1-u))$,
respectively.  Also $Y_{1,1}=U_{0,1}$ is uniformly distributed on $[0,1]$.
\end{proof}

\subsection{Uniqueness  of solutions of  \eqref{functional1}.}
Let    $\mu=\pa{\mu_\gl}_{\gl>0}$ and 
$\theta=\pa{\theta_\gl}_{\gl>0}$
be two solutions of \eqref{functional1} in $\mathcal M$.
Let  $Y=\pa{Y(\gl)}_{\gl>0}$  
and  $Z=\pa{Z(\gl)}_{\gl>0}$  denote  two measurable
processes representing respectively $\mu$ and 
$\theta$, in the sense of Remark \ref{remmesure} (i).
Without loss of generality, we can assume that $Y$ and $Z$ share the same
underlying  probabilistic space, and the same exponent $\alpha$. 
Then,  by definition of $\mathcal M$,
for $\Lambda>0$,
\[
d_\Lambda(Y,Z)
=
\sup_{(0,\Lambda)}\esp{\lambda^\alpha\vert Y(\lambda)-Z(\lambda)\vert} 
\]
is finite. Let $\delta$   
denote the infimum of $d_\Lambda(\hat Y,\hat Z)$
over all couples of representations $(\hat Y,\hat Z)$ of $\mu$ and $\theta$,
 lying on the same probabilistic space,
and assume that $\delta>0$. Let $(Y_0,Z_0)$ be
 such a couple of representations, 
satisfying furthermore
\[
d_\Lambda(Y_0,Z_0)
<
\delta\ \frac{3-\alpha}{2}.
\]
Consider a probabilistic space on  which are defined 
three \textit{independent} random variables
  $(Y_1,Z_1)$,  $(Y_2,Z_2)$ and $U$,
$(Y_1,Z_1)$ and  $(Y_2,Z_2)$ being two copies of $(Y_0,Z_0)$, 
$U$ being uniform  on  $(0,1)$. Finally,
 for every $\lambda>0$,  set 
\begin{align*}
\hat Y(\gl)
=
U^2 Y_1(\gl U) + (1-U)^2  Y_2(\gl(1-U))
+ \Theta(\gl,U),
\\
\hat Z(\gl)
=
U^2 Z_1(\gl U) + (1-U)^2  Z_2(\gl(1-U))
+ \Theta(\gl,U).
\end{align*}
Then $\hat Y$   and $\hat Z$  are   
representations of $\mu$ (resp. $\theta$) and satisfy 
Remark \ref{remmesure}  (ii). 
Moreover, we have, for $\lambda\in(0,\Lambda)$,
\begin{align*}
\esp{\lambda^\alpha\left|\hat Y(\lambda)-\hat Z(\lambda)\right|}
&=
\mathbb{E}\left[\lambda^\alpha\left|U^2 Y_1(\gl U) + (1-U)^2 Y_2(\gl(1-U))\right.\right.\\
&\hspace{2cm}
\left.\left.- U^2 Z_1(\gl U) - (1-U)^2 Z_2(\gl(1-U))\right|\right]
\\
&\leq 2\esp{\lambda^\alpha U^2\left|Y_1(\gl U)-Z_1(\gl U)\right|}
\\
&\leq 2\int_0^1 u^{2-\alpha}\ 
\esp{(\lambda u)^\alpha \left|Y_1(\gl u)-Z_1(\gl u)\right|}  du
\\
&<\delta,
\end{align*}
leading to a contradiction.

\subsection{Moments of $X(\gl)$.}
The aim of  this  Section is the computation  of moments of $X(\gl)$,  
 completing the proof of Theorem  \ref{p=1/n}.  
If one uses  directly \eqref{functional1}, the
computations of moments by induction
  are hardly tractable because all three terms
on the right of  \eqref{functional1} depend on $U$.  To circumvent
this problem,
  we consider  a new distributional identity
\begin{equation}
\label{functional2}
W(\gl)\egenloi
\xi(\gl)+UW(\gl U)+(1-U)\widetilde{W}(\gl(1-U)),
\end{equation}
in which
\begin{itemize}
\item $\xi(\gl)$ is as in Section \ref{S:results}; equivalently,
$\xi(\gl)=\sum_{x\in\Pi_1} x;
$
\item  $\xi(\gl)$ and $\pa{U,W(\gl U),\widetilde{W}(\gl(1-U))}$
are independent;
\item conditionally, given $U=u$,  $W(\gl U)$ and $\widetilde{W}(\gl(1-U))$
are independent and distributed as  $W(\gl u)$ and ${W}(\gl(1-u))$,
respectively.
\end{itemize}
The next Propositions establish relations between
$X(\gl)$ and  solutions of \eqref{functional2}, 
 eventually providing an algorithm for the computation of 
 moments of $X(\gl)$
(see \eqref{E:induction} and \eqref{E:induction2}).

\begin {prop}
\label{5to34}
The family
$\pa{Y(\gl)}_{\gl> 0}=\pa{\xi(\gl)+\gl X(\gl)}_{\gl> 0}$, in
which $\xi(\gl)$ and $X(\gl)$ are assumed independent,
  is a solution of \eqref{functional2}.
\end {prop}

\begin {prop}
\label{unique2}
The $n$-th moment  of $Y(\gl)$    is a polynomial of degree $n$
in the variable $\gl$.
\end {prop}

Before proving Propositions
\ref{5to34} and \ref{unique2}, we need a lemma.
\begin {lem}
\label{factsaboutW}
The $n$-th moment $g_n(\gl)=\esp{\xi(\gl)^n}$
is a polynomial of degree $n$ with nonnegative coefficients and for $n\geq1$,
$g_n(0)=0$.
\end {lem}
\begin{proof}
Owing to Campbell's Theorem \cite[p.28]{PPP},
we have
\begin{equation*}
\esp{ e^{s\xi(\gl)}} = \exp \pa{ \gl  \pa{  \esp{e^{sU}} -1}    }
=\exp \pa{ \gl  \pa{ \frac s{2!}+\frac {s^2}{3!}+\dots}    }.
\end{equation*}
Expanding the last expression gives the lemma.
\end{proof}

\begin{proof}[Proof of Proposition \ref{5to34}]

To show \eqref{functional2}, it is enough to show
\begin{align*}
\gl X(\gl)
&\egenloi
UY(\gl U)+(1-U)\widetilde{Y}(\gl(1-U))
\\
&\egenloi U\xi(\gl U)+\gl U^2X(\gl U)
 +(1-U)\widetilde\xi(\gl(1-U))+\gl (1-U)^2\widetilde{X}(\gl(1-U)),
\end{align*}
where, as usual, conditioned on $U=u$,
the terms on the right hand side are independent with the
right distributions.
This follows immediately from \eqref{functional1},
since
$$
\gl\Theta(\gl,u)
\egenloi
u\xi(\gl u) +(1-u)\widetilde\xi(\gl(1-u)).
$$
\end{proof}

\begin{proof}[Proof of Proposition \ref{unique2}]
Consider the sequence of integral equations
\begin{equation}
\label{E:induction}
P_0(\gl)=1,\qquad P_n(\gl)
=2\int_0^1u^nP_n(\gl u)du+\psi_n(\gl),
\quad n\ge 1,
\end{equation}
in which
\begin{equation}
\label{psi}
\psi_n(\gl)
=\sum_{r+k+\ell=n\atop{k<n,\ell<n}}\binom{n}{r,k,\ell}
g_r(\gl)
\int_0^1u^k(1-u)^{\ell}P_k(\gl u)P_{\ell}(\gl(1-u))\, du,
\end{equation}
where $g_r$ is the $r$-th moment of $\xi(\lambda)$.
Proposition \ref{unique2} is a consequence of the next lemma.
\end{proof}
\begin{lem}
The  induction formula \eqref{E:induction} and the initial condition
$P_1(0)=0$
 defines a unique sequence
 of  polynomials, $\pa{P_n(\gl)}_{n\ge 0}$.
Furthermore, $P_n$ has degree $n$,
and vanishes at $0$.
For $n\ge 1$,
the $n$-th  moment  $\esp{Y(\gl)^n}$  is  equal to $P_n(\gl)$.
\end{lem}

\begin{proof} Consider $n\ge 1$ and assume
 that the properties in the lemma hold for
$1\le m\le n-1$. Then, for $k$ and $\ell$ smaller than $n$,
and $r+k+\ell=n$, the expression
\[g_r(\gl)\int_0^1u^k
(1-u)^{\ell}P_k(\gl u)P_{\ell}(\gl(1-u))\, du\]
is a polynomial with degree $n$ and,
due to Lemma \ref{factsaboutW},
vanishes at 0.  Thus, in this case,   $\psi_n(\gl)$  is a
polynomial with degree $n$,
 vanishing at 0.
It is now easy to check that a polynomial $P_n(\gl)$
satisfies \eqref{E:induction} if and only if, for
$(n,k)\neq (1,0)$,
\begin{equation}
\label{fn2}
\cro{\gl^k}P_n=\frac{n+k+1}{n+k-1}\ \cro{\gl^k}\psi_n.
\end{equation}
Also, by the induction assumptions,
\begin{align*}
\esp{Y(\gl)^n}
&=
\esp{\big(\xi(\gl)+UY(U\gl)+(1-U)\widetilde{Y}((1-U)\gl)\big)^n}\\
&=\sum_{r+k+\ell=n}
\binom{n}{r,k,\ell}g_r(\gl)
\esp{U^k(1-U)^{\ell}Y(U\gl)^k\widetilde{Y}((1-U)\gl)^{\ell}}\\
&=2\esp{U^nY(\gl U)^n}+\psi_n(\gl).
\end{align*}
Note that $\psi_n(\gl)\ge 0$ for $\gl\ge 0$.
By Remark \ref{remmesure},
 $\gl\to f_n(\gl)=\esp{Y(\gl)^n}$ is nonnegative and
measurable. Thus, for $\gl>0$, we can
rewrite the previous equation:
\begin{align*}
f_n(\gl)
&=
2\int_0^1u^nf_n(\gl u)du+\psi_n(\gl)\\
&=2 \gl^{-n-1}\int_0^\gl v^nf_n(v)dv+\psi_n(\gl).
\end{align*}
Since $f_n(\gl)$  is  assumed to be finite and $\psi_n(\gl)\ge 0$,
  the  integral on the right hand side is convergent, and thus it is a continuous function
 of $\gl$.
As a consequence  $f_n$ belongs to  $C^\infty(0,+\infty)$, and is a solution 
 on $(0,+\infty)$ of the following
differential equation:
$$\gl f_n'(\gl)+(n-1)f_n(\gl)=(n+1)\psi_n(\gl)+\gl\psi_n'(\gl).$$
by Proposition \ref{ilexiste1} and Lemma \ref{factsaboutW},
  $\gl^{n-1}f_n(\gl)\to0$ as $\gl\to0$, but the general solution of the
differential equation is
$P_n(\gl)+C\, \gl^{-n+1}$. Thus $f_n=P_n$  on $(0,+\infty)$.
\end{proof}
As a consequence of these results, we deduce that:
\begin{prop}
\label{polynomialprop}
The function $\gl\longrightarrow\gl^n\esp{X(\gl)^n}$ is a polynomial of degree $n$ that vanishes 
at $0$.
\end{prop}

\begin{proof}
Since $Y(\gl)=\xi(\gl)+\gl X(\gl)$, with independent summands,  we obtain
\begin{equation}
\label{E:induction2}
P_m(\gl)
=
\esp{Y(\gl)^m}
=
\sum_{0\le k\le m}{m\choose k}\esp{\xi(\gl)^{m-k}}\gl^k\esp{ X(\gl)^k}.
\end{equation}
The result follows by induction.   
\end{proof}

\subsection{Computation of the first moments}
The moments of $Y(\gl)$, and thus of $X(\gl)$, can be computed up to arbitrary order
with the help of \eqref{E:induction} and  \eqref{E:induction2}.
For the first two moments, the calculations run as follows.
Expanding
\[\exp \pa{ \gl  \pa{ \frac s{2!}+\frac {s^2}{3!}+\dots}    },\]
in the proof of   Lemma \ref{factsaboutW} we obtain
\begin{lem}\label{L:xi}
$$g_1(\gl)=\esp{\xi(\gl)}=\frac12\gl
\qquad\text{and}\qquad
g_2(\gl)=\esp{\xi(\gl)^2}=\frac13\gl+\frac14\gl^2.$$
\end{lem}

\begin{prop}
$$\gl\esp{X(\gl)}=\esp{\Xi(\gl)}=\gl
\qquad\text{and}\qquad
\gl^2\var{X(\gl)}=\var{\Xi(\gl)}=\frac13\gl+\frac{1}{12}\gl^2.$$
\end{prop}

\begin{proof}
Taking $n=1$ in \eqref{psi} and \eqref{fn2}, we find, using Lemma \ref{L:xi},
\begin{align*}
\psi_1(\gl)&=g_1(\gl) = \tfrac12\gl,
\\
P_1(\gl)&= \tfrac31\cdot\tfrac12\gl=\tfrac32 \gl.
\end{align*}
Taking $n=2$, we similarly find
\begin{align*}
\psi_2(\gl)
&=g_2(\gl)
+ 2\cdot 2g_1(\gl) \int_0^1 uP_1(\gl u)\,du
+2\int_0^1 u(1-u)P_1(\gl u)P_1(\gl(1-u))\,du
\\
&=
\tfrac13\gl+\tfrac75\gl^2,
\\P_2(\gl)&=
\tfrac42\cdot\tfrac13\gl+\tfrac53\cdot\tfrac75\gl^2
=\tfrac23\gl+\tfrac73\gl^2.
\end{align*}
Since $Y(\gl)=\xi(\gl)+\gl X(\gl)$, with independent summands,
\[P_1(\gl)=\esp{Y(\gl)}=\esp{\xi(\gl)}+\gl\esp{ X(\gl)}, \]
which by Lemma \ref{L:xi}
yields $\gl\esp{X(\gl)}=\gl$. Similarly,
\begin{align*}
\gl^2\esp{X(\gl)^2} = P_2(\gl)-\esp{\xi(\gl)^2}-2\esp{\xi(\gl)}\esp{\gl X(\gl)}
=\tfrac13\gl+\tfrac{13}{12}\gl^2,
\end{align*}
which yields the variance formula.
\end{proof}

The formulas for mean and variance of $X(\gl)$ can also be obtained
directly from
\eqref{samuel} and Lemma \ref{L:xi}; we leave this as an exercise.

\section{Concluding remarks}

We have presented a probabilistic analysis of Quicksort
when some comparisons
can err. Analysing other sorting algorithms such as merge sort,
insertion sort or selection
is even more intricate. They do not fit into the model presented in
this paper and further
more involved probabilistic models/arguments are required.
We conjecture  that the same normalization holds for the number of inversions
in the \textit{output of merge sort}
for $n=2^m\rightarrow +\infty$, $p=\gl/n$, and that the
   limit law $\widehat X(\gl)$ satisfies
\[\esp{\widehat X(\gl)}=\sum_{k\ge 0}\ \frac{2^k}{(2^k+2)(2^k+3)}=
0.454674373\dots<\esp{X(\gl)}.\]

\section{Acknowledgements}

A more rigorous formulation  and a shorter proof    of  
Theorem \ref{thm:unique} are born from discussions with 
 Uwe R\"osler. Also, we thank 
two anonymous  referees,
 whose careful reading led to substantial improvements.

\end{document}